\documentclass[11pt]{amsart}
\usepackage{geometry}                % See geometry.pdf to learn the layout options. There are lots.
\usepackage[parfill]{parskip}    % Activate to begin paragraphs with an empty line rather than an indent
\usepackage{graphicx}
\usepackage{amssymb}
\usepackage{epstopdf}
\usepackage{xcolor}
\usepackage{amsthm}
\usepackage{amsmath}
\usepackage{hyperref}
\usepackage{todonotes}
\usepackage{comment}
\usepackage{tikz}
\newenvironment{proof1}[1]{\begin{trivlist} \item[] {\em Proof of #1:}}{\newline \textcolor{white}{.}\hfill $\Box$
                      \end{trivlist}}
\theoremstyle{plain}

\newtheorem{thm}{Theorem}[section]
\newtheorem{lem}{Lemma}[section]
\newtheorem{prop}{Proposition}[section]
\newtheorem{cor}{Corollary}[section]

\newcommand{\pa}{\partial}

\newcommand{\R}{\mathbb{R}}

\newcommand{\Hone}{H^1_{\gamma}(\Omega)}
\newcommand{\Lpa}{L^2(\pa\Omega)}
\newcommand{\Hplus}{H^{1/2}_{\gamma}(\pa\Omega)}
\newcommand{\Hminus}{H^{-1/2}_{\gamma}(\pa\Omega)}
\newcommand{\norm}[1]{\left\lVert#1\right\rVert}

\newcommand{\XM}{{X_{\gamma}(M)}}
\newcommand{\HN}[1]{{\mathbb{H}^{#1}_{\gamma}}}
 \newcommand{\PsiN}[1]{{{\Psi}_{\gamma}^{#1}}}

\author{Thomas Beck}
\email{tbeck7@fordham.edu}
\address{Mathematics Department, Fordham University, Bronx, NY 10458}
\author{Yaiza Canzani}
\email{canzani@email.unc.edu}
\address{Department of Mathematics, University of North Carolina, Chapel Hill, NC 27514}
\author{Jeremy L. Marzuola}
\email{marzuola@math.unc.edu}
\address{Department of Mathematics, University of North Carolina, Chapel Hill, NC 27514}

\vspace{-0.3in}
\title[Non-local exchange operator estimates]{On non-local exchange and scattering operators in domain decomposition methods}
\date{\vspace{-0.4in}\today}
\begin{document}

\begin{abstract}
We study non-local exchange and scattering operators arising in domain decomposition algorithms for solving elliptic problems on domains in $\mathbb{R}^2$. Motivated by recent formulations of the Optimized Schwarz Method introduced by Claeys, we rigorously analyze the behavior of a family of non-local exchange operators $\Pi_\gamma$, defined in terms of boundary integral operators associated to the fundamental solution for $-\Delta + \gamma^{-2}$, with $\gamma > 0$. Our first main result establishes precise estimates comparing $\Pi_\gamma$ to its local counterpart $\Pi_0$ as $\gamma \to 0$, providing a quantitative bridge between the classical and non-local formulations of the Optimized Schwarz Method. In addition, we investigate the corresponding scattering operators, proving norm estimates that relate them to their classical analogues  through a detailed analysis of the associated Dirichlet-to-Neumann operators. Our results clarify the relationship between classical and non-local formulations of domain decomposition methods and yield new insights that are essential for the analysis of these algorithms, particularly in the presence of cross points and for domains with curvilinear polygonal boundaries.
\end{abstract}

\maketitle
%\vspace{-0.2in}

\section{Introduction}

Domain decomposition methods have long been recognized as efficient and effective approaches for solving elliptic equations. Among these, the Optimized Schwarz Method has emerged as a popular and versatile technique. Introduced by Lions \cite{lions1990schwarz} and further developed by Després \cite{despres1990,despres1991domain,despres1991-2,despres1993}, this iterative procedure solves subdomain problems whose solutions converge to the global solution on the entire domain. The method has been successfully applied in various settings, including the numerical approximation of the Helmholtz equation. For additional developments and applications, see \cite{collino2000domain,gander2001optimized,gander2006optimized} and the references therein.

Suppose one wishes to solve $\Delta u=f$ on a set $\Omega \subset \R^2$.  If we decompose $\Omega$ as the union of two sets $\Omega_0, \Omega_1$ such that their intersection $\Gamma=\Omega_0 \cap \Omega_1$ is a separating curve, then the Optimized Schwarz Method is an iterative scheme for finding a solution $u$ in $\Omega$ given by $u_j$ on $\Omega_j$, $j=0,1$, where 
$u_j = \lim_{n\to \infty} u_j^n$ and the $u_j^n$ solve the system
 \begin{align*}
     &  \Delta u_j^{n+1} = f \ \ \text{in} \ \ \Omega_j, \qquad j=0,1,\\
     &  F_0(u_0^{n+1},\partial_\nu u_0^{n+1}) = F^*_0(u_1^{n},\partial_\nu u_1^{n}) \ \ \text{on} \ \ \Gamma, \\
  %   &  \Delta u_2^{n+1} = f \ \ \text{in} \ \ \Omega_2, \\
     &   F_1(u_1^{n+1},\partial_\nu u_1^{n+1}) = F^*_1(u_0^{n},\partial_\nu u_0^{n}) \ \ \text{on} \ \ \Gamma.
 \end{align*}
Here, each $F_j$ is a linear operator on the Dirichlet and Neumann traces on $\Gamma$, and the method is initialized by setting, for example,   $u_0^0 = u_1^0 = 0$.
A standard approach is to let $F(u,v) = v + S(u)$, $F^* (u,v) = -v + S^* (u)$, where $S$ is a linear operator.  When $S(u)=u$, this is the Robin-type condition as originally proposed by Lions, \cite{lions1990schwarz}.  Other linear operators constructed on the artificial boundary $\Gamma$ are possible, and the choice of operators $F,F^*$ that one uses on the traces can determine the efficacy of the method and its rate of convergence.  See for instance \cite{claeys2022robust},  Section 5, for a variety of examples. 
When the domains overlap, unlike in our setup, it is common to choose $F(u,v) = u, F^*(u,v) = u$ as a Dirichlet matching condition, see \cite{gander2006optimized}.

 Using interior elliptic regularity of the solution, the solutions should satisfy $u_0 = u_1$ on $\Gamma$, as well as $\partial_{\nu} u_0 = - \partial_{\nu} u_1$ on $\Gamma$. To ensure the latter, many formulations of the Optimized Schwarz Method hinge upon using an exchange operator, $\Pi_0$, defined on the shared boundary $\Gamma$, which switches the Neumann traces on the boundary of the domains. The switching is done up to a sign, coming from the normal vectors pointing in opposite directions across the shared boundary. Note that this exchange operator acts locally on the Neumann traces. See \eqref{eqn:Pi0} for the precise definition in the case where $\Omega$ is decomposed into two subdomains. For example, such a local exchange operator is used in the original algorithms of Despr\'es, \cite{despres1990}, \cite{despres1991-2}, as well as in \cite{gander2002}, and the general presentation of non-overlapping domain decomposition methods from \cite{collino2000domain}.

While the Optimized Schwarz Method is well known to converge numerically for discretized problems across a wide range of transmission coefficients (see \cite{lions1988schwarz, lions1990schwarz, despres1991domain}), its rate of convergence is highly sensitive to the mesh size and structure \cite{bendali2006non}. In practice, achieving stability becomes particularly challenging near cross-points, where multiple subdomains intersect, and careful treatment is required to maintain accuracy \cite{despres2021corners}. Notably, Gander and Kwok \cite{gander2013applicability} demonstrated that the convergence results established by Lions at the continuous level \cite{lions1990schwarz}, specifically in $H^1(\Omega_j)$, do not necessarily extend to the discretized setting, where finite element spaces are employed.

To overcome these difficulties and ensure both stability and convergence in more complex decompositions, Claeys \cite{Cl} and Claeys-Parolin \cite{claeys2022robust} proposed a powerful alternative formulation of the Optimized Schwarz Method for the Helmholtz problem. Their approach is specifically designed to handle decompositions with cross-points, where several subdomains meet (see Figure 1 for a two-dimensional example recreated from \cite{Cl}). In Section 7 of \cite{Cl}, their method is shown to be strongly coercive at the continuous level. In Section 11 of \cite{claeys2022robust}, they present a numerical implementation that demonstrates geometric convergence of the discretized method, with rates that are uniform with respect to the mesh size.

This robust formulation has since been extended to a variety of challenging settings, including non-self-adjoint impedance operators \cite{claeys2023}, time-harmonic Maxwell equations \cite{ccp2022}, and the Helmholtz equation in bounded cavities \cite{claeys2023-2}. In addition, recent developments have significantly reduced the computational cost of applying these non-local exchange operators \cite{claeys2024}, making the method increasingly practical for large-scale computations. As shown in \cite{claeys2022robust}, when the domain partitioning avoids cross-points and an appropriate impedance operator is chosen, this new formulation reduces to the classical Optimized Schwarz Method of Despr\'es \cite{despres1991domain}, offering a unified framework.

 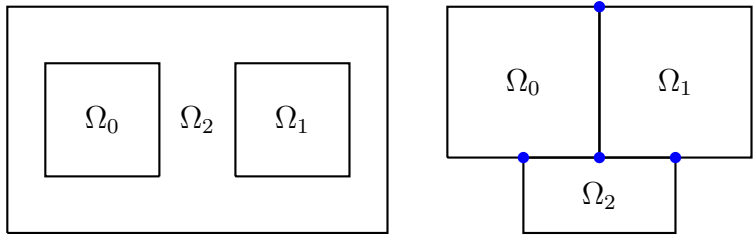
\begin{figure}[h]
    \centering
\begin{tikzpicture}
\draw[thick] (-5/2,-3/2) -- (-5/2,3/2) -- (5/2,3/2) -- (5/2,-3/2) -- (-5/2,-3/2);
\draw[thick] (1/2,-3/4) -- (1/2,3/4) -- (2,3/4) -- (2,-3/4) -- (1/2,-3/4);
\draw[thick] (-1/2,-3/4) -- (-1/2,3/4) -- (-2,3/4) -- (-2,-3/4) -- (-1/2,-3/4);
\node[font=\large] at (-5/4,0) {$\Omega_0$};
\node[font=\large] at (5/4,0) {$\Omega_1$};
\node[font=\large] at (0,0) {$\Omega_2$};
%\draw[<->] (0,-0.2)--(0,0.2);
%\node at (0,0.5) {$w$};
\end{tikzpicture}
\hspace{.5cm}
\begin{tikzpicture}
%\draw[thick] (-1,0.2) -- (1,0.2);
%\draw[thick] (-1,-0.2) -- (1,-0.2);
\draw[thick] (-2,-1) -- (-2,1) -- (0,1) -- (0,-1) -- (-2,-1);
\draw[thick] (0,-1) -- (0,1) -- (2,1) -- (2,-1) -- (0,-1);
\draw[thick] (-1,-1) -- (1,-1) -- (1,-2) -- (-1,-2) -- (-1,-1);
\node[font=\large] at (-1,0) {$\Omega_0$};
\node[font=\large] at (1,0) {$\Omega_1$};
\node[font=\large] at (0,-3/2) {$\Omega_2$};
\node [circle, fill=blue, minimum size = 0.15cm, inner sep=0pt] at (-1,-1) {};
\node [circle, fill=blue, minimum size = 0.15cm, inner sep=0pt] at (1,-1) {};
\node [circle, fill=blue, minimum size = 0.15cm, inner sep=0pt] at (0,-1) {};
\node [circle, fill=blue, minimum size = 0.15cm, inner sep=0pt] at (0,1) {};
%\draw[<->] (0,-0.2)--(0,0.2);
%\node at (0,0.5) {$w$};
\end{tikzpicture}
    \caption{Examples from \cite{Cl} of a sub-domain partition without cross-points (left) and with cross-points denoted with blue dots (right). }
    \label{fig:domdec}
\end{figure}

As for the original Optimized Schwarz Method, information needs to be exchanged between subdomains to, in particular, impose the condition $\partial_{\nu}u_0 = -\partial_{\nu}u_1$ across $\Gamma$. In \cite{Cl}, Claeys does this via the introduction of an \emph{exchange operator}, $\Pi_\gamma$, indexed by a parameter $\gamma > 0$. See \eqref{defn:exchange} for a precise definition in the case where $\Omega$ is decomposed into two subdomains. Unlike $\Pi
_0$, $\Pi_\gamma$ is an operator which acts non-locally on the Neumann traces, and is
 defined in terms of boundary integral operators that are built by integrating the fundamental solution $-\Delta+\gamma^{-2}$ on the boundary of each subdomain. Unlike for the local exchange operator $\Pi_0$, the operators $\Pi_{\gamma}$ couple the Neumann traces of even those subdomains that do not have any shared boundary.  It is in this manner that we refer to the methods here as non-local, though we acknowledge that the Neumann traces are themselves non-local in the more standard sense of the term.

One of the key advantages of Claeys’s formulation is its robustness: the method remains coercive even in the presence of cross-points in the domain decomposition. However, this benefit comes at a cost. The non-local nature of the exchange operator $\Pi_\gamma$ introduces significant challenges in analyzing and implementing the associated iterative scheme. This stands in contrast to the Optimized Schwarz Method, where, in the absence of cross-points, the iterative scheme based on the simpler, local exchange operator $\Pi_0$ is straightforward to set up and execute.

A primary goal of this work is to rigorously establish that $\Pi_\gamma$ converges to $\Pi_0$ as $\gamma \to 0$, in an appropriate functional sense. Our first main result, Theorem \ref{thm:exchange} below,  builds a precise quantitative bridge between the two operators, offering new insight into their relationship. This asymptotic behavior was previously conjectured by Claeys in private communication, and we confirm and extend his intuition through our analysis.

The formulations of the Optimized Schwarz Method using the exchange operator $\Pi_0$ also involve \emph{scattering operators}, $S^j_0$, which map ingoing to outgoing Robin traces of solutions to a Helmholtz problem in each subdomain $\Omega_j$. These operators act locally on the Dirichlet and Neumann traces of each subdomain. Scattering operators in domain decomposition methods for Helmholtz problems have for example appeared in \cite{despres1991-2,lecouvez2015iterative}. In \cite{Cl}, Claeys introduces a counterpart to these operators, $S^j_\gamma$, again indexed by the parameter $\gamma > 0$, and uses it in a crucial manner when reformulating the Helmholtz scattering problem as a well posed problem on the boundaries of the subdomains.  These scattering operators, defined precisely in \eqref{eqn:scattering1}, are expressed in terms of the Dirichlet-to-Neumann operator, $T^j_\gamma$, of $\Omega_j$, associated to the operator $-\Delta + \gamma^{-2}$. In particular, unlike $S_0^j$, the $S^j_\gamma$ act non-locally on the  Dirichlet and Neumann traces of the subdomains.
 Our second main result, Theorem \ref{thm:DtN}, provides estimates on $T^j_\gamma - Id$ for general curvilinear polygonal domains $\Omega_j$. These estimates, in addition to their intrinsic interest in the theory of Dirichlet-to-Neumann operators on non-smooth domains, lead directly to a norm estimate relating $S^j_\gamma$ and $S^j_0$ in the $\gamma \to 0$ limit (see Corollary \ref{cor:DtN}). This establishes a quantitative relationship between the classical formulations of the Optimized Schwarz Method and the new formulation introduced in \cite{Cl}.

In Section \ref{sec:Claeys}, we describe in detail the formulation of the method introduced by Claeys in \cite{Cl}. This includes the definitions of the non-local exchange and scattering operators for a general multi-domain decomposition, potentially including cross points. Although Theorems \ref{thm:exchange} and \ref{thm:DtN} are stated for the case of a decomposition into two subdomains, we show in Section \ref{sec:Claeys} how these results extend to the general multi-domain setting. In particular, our estimates provide quantitative bounds between the non-local operators and their local counterparts in the $\gamma \to 0$ limit, clarifying the connection between Claeys' formulation and the classical Optimized Schwarz Method.

\subsection{The exchange operator} \label{subsec:exchange}
 Let $\Omega_1$ be a bounded planar domain, with boundary given by a smooth, curvilinear polygon, and set $\Omega_0 = \R^2\backslash \bar{\Omega}_1$. We let $\Gamma = \pa\Omega_0= \pa\Omega_1$ correspond to the shared boundary, and for each $\gamma>0$ define the function spaces 
\begin{align*}
\HN{-1/2}(\Gamma) := H^{-1/2}_\gamma(\pa\Omega_0) \times  H^{-1/2}_\gamma(\pa\Omega_1), \qquad \HN0(\Gamma):= L^2(\pa\Omega_0) \times  L^2(\pa\Omega_1).
\end{align*}
The $H^{-1/2}_{\gamma}(\pa\Omega_j)$ spaces are defined by duality: $H_{\gamma}^1(\Omega_j)$ is given by $H^1(\Omega_j)$ functions with norm
    \begin{align*}
        \norm{F}_{H^1_{\gamma}(\Omega_j)} := \norm{\nabla F}_{L^2(\Omega_j)}+\gamma^{-2}\norm{F}_{L^2(\Omega_j)}.
    \end{align*}
The $H^{1/2}_{\gamma}(\pa\Omega_j)$ space is then the restriction of $F\in H_{\gamma}^1(\Omega_j)$ to $\pa\Omega_j$, with norm
    \begin{align*}
        \norm{f}_{H^{1/2}_{\gamma}(\pa\Omega_j)} := \inf\{ \norm{F}_{H^1_{\gamma}(\Omega_j)}: F|_{\pa\Omega}=f\}.
    \end{align*}
    The $H^{-1/2}_{\gamma}(\pa\Omega_j)$ space is the dual of $H^{1/2}_{\gamma}(\pa\Omega_j)$ and $\HN{s}(\Gamma)$  is defined by interpolation for $-\tfrac{1}{2} < s < 0$. Then, for $\mathcal{G}_\gamma$ a fundamental solution for $-\Delta+\gamma^{-2}$, and $\phi = (\phi_0,\phi_1)\in \HN{-1/2}(\Gamma)$, let
\begin{align} \label{eqn:Psi1}
    \PsiN{}(\phi)(x)% & := \PsiN{0}(\phi_0)(x) + \PsiN{1}(\phi_1)(x) \\ 
    & := \int_{\pa\Omega_0}\mathcal{G}_\gamma(x-y)\phi_0(y)\,d\sigma(y) + \int_{\pa\Omega_1}\mathcal{G}_\gamma(x-y)\phi_1(y)\,d\sigma(y).
\end{align}
Note that $\mathcal{G}_{\gamma}(x) =K_0(|x|/\gamma)$ where $K_0$ is the order $0$ modified Bessel function of the second kind. Here $d\sigma$ is the line measure on $\pa\Omega_0 = \pa\Omega_1$.

The \emph{exchange operator},  $\Pi_\gamma$, acts on $\phi=(\phi_0,\phi_1)\in\HN{-1/2}(\Gamma)  \times \HN{-1/2}(\Gamma) $ and is defined by
\begin{equation}\label{defn:exchange}
    \Pi_\gamma \begin{pmatrix} \phi_0\\\phi_1 \end{pmatrix}
    := \begin{pmatrix} \phi_0\\\phi_1 \end{pmatrix} 
    -2\begin{pmatrix} \pa_{n_0}^{int}\,\PsiN{}(\phi)\\\pa_{n_1}^{int}\,\PsiN{}(\phi) \end{pmatrix},
\end{equation}
where $\PsiN{}$ is defined in \eqref{eqn:Psi1}. Here, $n_j$ is the outward pointing normal to $\pa\Omega_j$, and the superscript $int$ refers to these derivatives being evaluated from the interior of $\Omega_j$ for $j=0,1$.

Theorem \ref{thm:exchange} approximates $\Pi_\gamma$ by its local counterpart $\Pi_0$,
\begin{align}\label{eqn:Pi0}
    \Pi_0\begin{pmatrix} \phi_0\\\phi_1 \end{pmatrix} =- \begin{pmatrix} \phi_1\\\phi_0 \end{pmatrix}.
\end{align}

As we will show in Section \ref{sec:exchange} below, $\Pi_\gamma$ and $\Pi_0$ are linked via a single layer potential.

As discussed above, $\Pi_\gamma$ is precisely the non-local exchange operator introduced by Claeys in \cite{Cl} in the case of two subdomain decomposition, while $\Pi_0$ is the local exchange operator from classical formulations of the Optimized Schwarz Method. Our first main theorem gives a relationship between these non-local and local exchange operators.  When $\Gamma = \pa\Omega_0=\pa\Omega_1$ is smooth, this estimate will hold on the whole of $\HN{s}(\Gamma)$, while if $\Gamma$ contains a vertex, then the estimate will be restricted to those functions in $\HN{s}(\Gamma)$ with a bound on their Dirichlet energy on $\Gamma$, coming from the following definition:

For each  $\gamma>0$ and $0< M\leq \gamma^{-1}$, we define the set 
\begin{align} \label{eqn:XM}
    \XM := \{f\in H^{1/2}_{\gamma}(\Gamma)\cap H^1(\Gamma): \norm{\nabla f}_{L^2(\Gamma)} \leq M \norm{f}_{L^2(\Gamma)}\}.
    \end{align}
Similar spaces and related regularity restrictions appear in the author's previous work \cite{beck2022quantitative} on impedance-to-impedance bounds for Poincar\'e-Steklov hierarchical numerical schemes used in solving Helmholtz problems.  

        The set $\XM$ contains a subspace of $H^{1/2}_{\gamma}(\Gamma)$ of dimension comparable to $M$. Given $f\in H^1(\Gamma)$, there exists $c>0$ such that for all $0<\gamma<c$, $f$ will be in $\XM$ for $M$ sufficiently large.       

\begin{thm} \label{thm:exchange}
Let $\Omega_1$ be a domain in $\R^2$. There exist constants $c$, $C$ such that the following hold.
\begin{itemize}
    \item[1.] If $\Omega_1$ has smooth boundary, then for all $0<\gamma<c$, $-\tfrac{1}{2}\leq s \leq 0$, and $\phi \in \HN{s}(\Gamma)$,
\begin{align*}
    \norm{(\Pi_\gamma-\Pi_0)\phi}_{\HN{s}(\Gamma)} \leq C\gamma\norm{\phi}_{\HN{s}(\Gamma)}.
\end{align*}
\item[2.] If $\Omega_1$ is a smooth curvilinear polygon, then for each $0<a<1$, there exists a constant $C_a>0$ such that for all $0<\gamma<c$, $0<M<C_a^{-1}\gamma^{-a}$,  $-\tfrac{1}{2}\leq s \leq 0$, and $\phi = (\phi_0,\phi_1) \in \HN{s}(\Gamma)$ with $\phi_0$, $\phi_1\in \XM$,
\begin{align*}
    \norm{(\Pi_\gamma-\Pi_0)\phi}_{\HN{s}(\Gamma)} \leq C(1+M^{1/2})\gamma^{1/2}\norm{\phi}_{\HN{s}(\Gamma)}.
\end{align*}
\end{itemize}
\end{thm}

    Notably, the estimate in item 2 of Theorem \ref{thm:exchange} cannot be extended to the whole of $\HN{s}(\Gamma)$. We will show in Section \ref{sec:exchange} that there exists $c>0$, independent of $\gamma$, such that for each $0<\gamma<c$,  $-\tfrac{1}{2}\leq s \leq 0$, there is $\phi_\gamma\in \HN{s}(\Gamma)$ with
    \begin{align}\label{rem:exchange}
        \norm{(\Pi_\gamma-\Pi_0)\phi_\gamma}_{\HN{s}(\Gamma)} \geq c\norm{\phi_\gamma}_{\HN{s}(\Gamma)}.
    \end{align}
    If $\Gamma$ contains a vertex, then since $\Omega_0 = \mathbb{R}^2\backslash\bar{\Omega}_1$, at least one of $\Omega_0$ and $\Omega_1$ will have a convex vertex. It is the presence of this convex vertex that causes this restriction to $\phi_0$, $\phi_1\in \XM$. We will prove Theorem \ref{thm:exchange} in Section \ref{sec:exchange} by first writing $\Pi_\gamma$ in terms of an integral operator and then establishing bounds for its integral kernel.

\subsection{The scattering operator} \label{subsec:DtN}
In this section we study the local and non-local scattering operators appearing in Optimized Schwarz Methods.  In order to do so, we will first establish an estimate on a shifted Dirichlet-to-Neumann (DtN) operator of a, possibly unbounded, domain $\Omega\subset \R^2$, with boundary given by a piecewise smooth curvilinear polygon. 

For $\gamma>0$, we study the DtN operator $T_{\gamma}:H^{1/2}_{\gamma}(\pa\Omega)\to H^{-1/2}_{\gamma}(\pa\Omega)$ associated to $-\Delta+\gamma^{-2}$. That is, $$T_{\gamma}f = \gamma \pa_{\nu}F,$$ where
    $$\begin{cases}
        (-\Delta+\gamma^{-2})F  = 0 & \text{ in }\Omega, \\
        F|_{\pa\Omega}   = f & \text{ on }\pa\Omega.
    \end{cases}$$
As for the analysis of the exchange operator $\Pi_\gamma$, 
    the $H^{-1/2}_{\gamma}(\pa\Omega)$ space is the dual of $H^{1/2}_{\gamma}(\pa\Omega)$, defined via $H^1_\gamma(\Omega)$. The DtN operator, $T_\gamma$, can be used to define the $H^{1/2}_{\gamma}(\pa\Omega)$ and $H^{-1/2}_{\gamma}(\pa\Omega)$ scalar products
    \begin{align*}
        (f,g)_{H^{1/2}_{\gamma}(\pa\Omega)}  = \langle \gamma^{-1}T_{\gamma}f,\bar{g}\rangle_{\pa\Omega}, \qquad
        (p,q)_{H^{-1/2}_{\gamma}(\pa\Omega)}  = \langle p,\gamma T^{-1}_{\gamma}\bar{q}\rangle_{\pa\Omega}
    \end{align*}
    with the pairing $\langle f,g\rangle_{\pa\Omega} = \int_{\pa\Omega} fg\,d\sigma$. The latter can also be used to define the $H^{-1/2}_{\gamma}(\pa\Omega)$-norm.

Our main result concerning $T_\gamma$, quantifies how closely $T_\gamma$ approximates the identity operator when $\gamma > 0$ is small.  However, since $T_\gamma$ has a sequence of eigenvalues diverging to infinity for each fixed $\gamma$, such an estimate cannot hold uniformly on the entire space $H^{1/2}_\gamma(\partial \Omega)$.

Moreover, when the boundary $\partial \Omega$ includes convex vertices, we demonstrate that $T_\gamma$ possesses a finite set of eigenvalues, determined by the corresponding convex angles, which remain uniformly bounded away from $1$. Importantly, we also show that the eigenfunctions associated with these eigenvalues exhibit large Dirichlet energy.

As a result, and in analogy with Theorem \ref{thm:exchange}, we obtain a precise estimate for $T_\gamma$, with domain restricted to those functions in $H^{1/2}_{\gamma}(\pa\Omega)$ with a bound on their Dirichlet energy on $\pa\Omega$, coming from the function space $\XM$ in \eqref{eqn:XM}, with $\Gamma = \pa\Omega$.

    \begin{thm} \label{thm:DtN}
    Let $\Omega$ be a a smooth curvilinear polygon in $\R^2$. There exists a constant $c>0$ such that for all $\gamma \in (0,c)$, the following holds.
    \begin{enumerate}
        \item[1.] If $\Omega$ has no convex vertices, then there exists a constant $C$ such that for all $0<M<C^{-1}\gamma^{-1}$ and $h\in \XM$
 \begin{align*}
   \norm{(T_{\gamma}-Id)h}_{\Hminus} \leq C(\gamma^{1/2} + M^{1/2}\gamma^{1/2})\norm{h}_{\Hminus}.
    \end{align*}
        \item[2.] If $\Omega$ has at least one convex vertex, then for each $0<a<1$ there exists a constant $C_a$ such that for all $0<M <C_a^{-1}\gamma^{-a}$ and $h\in \XM$
    \begin{align*}
        \norm{(T_{\gamma}-Id)h}_{\Hminus} \leq C_a(\gamma^{1/3} + M^{1/2}\gamma^{a/2})\norm{h}_{\Hminus}.
    \end{align*}
    \item[3.]
           If $\Omega$ has smooth boundary, then there exists a constant $C$ such that for all $0<M<C^{-1}\gamma^{-1}$ and $h\in \XM$
         \begin{align*}
        \norm{(T_{\gamma}-Id)h}_{\Hminus} \leq C(\gamma^{1/2} + M\gamma)\norm{h}_{\Hminus}.
    \end{align*}
    \end{enumerate}
        
    \end{thm}

    We prove Theorem \ref{thm:DtN} in Section \ref{sec:DtN}, using estimates on the Robin eigenvalues for large parameters, and an upper bound on a Steklov Rayleigh quotient for functions in $\XM$. In the case of no convex vertices, these eigenvalue estimates will come from \cite{pank2013}, while we will use results from  \cite{Kh} when $\Omega$ has convex vertices. These works, for example, extended previous Robin eigenvalue estimates in terms of the maximum curvature of a smooth boundary, \cite{emp2014}, \cite{pank2015}.

Next, we explain how to apply Theorem \ref{thm:DtN} to the study of the non-local scattering operator introduced by Claeys in \cite{Cl}. For a function $u\in H^1_\gamma(\Omega)$, denoting $\tau_N(u)$ and $\tau_D(u)$ to be the Neumann and Dirichlet traces of $u$ on $\pa\Omega$, local and non-local ingoing and outgoing Robin traces are defined by
\begin{align*}
    \tau_{\pm,0}(u) = \tau_N(u) \pm i\tau_D(u), \qquad \tau_{\pm,\gamma}(u) = \tau_N(u) \pm i T_\gamma(\tau_D(u)).
\end{align*}
This then gives the definition of the scattering operators $S_0$ and $S_\gamma$ via
\begin{align} \label{eqn:scattering1}
    S_0(\tau_{-,0}(u)) = \tau_{+,0}(u), \qquad   S_\gamma(\tau_{-,\gamma}(u)) = \tau_{+,\gamma}(u).
\end{align}
Here, $S_0$ acts locally on the Dirichlet and Neumann traces of $u$, while $S_\gamma$ is the non-local version used in \cite{Cl}. Since
\begin{align*}
    S_0(\tau_{-,0}(u)) - S_\gamma(\tau_{-,\gamma}(u)) = \tau_{+,0}(u) - \tau_{+,\gamma}(u) = i\tau_{D}(u) - iT_\gamma(\tau_D(u)),
\end{align*}
Theorem \ref{thm:DtN} gives the following immediate relationship between $S_0$ and $S_\gamma$. 
\begin{cor} \label{cor:DtN}
   Let $\Omega$ be a smooth curvilinear polygon in $\R^2$. There exists a constant $c>0$ such that for all $\gamma \in (0,c)$, the following holds.
    \begin{enumerate}
         \item[1.] If $\Omega$ has no convex vertices,  then there exists a constant $C$ such that for all $0<M<C^{-1}\gamma^{-1}$, and all $u\in H^1_\gamma(\Omega)$ such that $\tau_D(u)$ is in $\XM$,
 \begin{align*}
   \norm{(S_0(\tau_{-,0}(u))-S_\gamma(\tau_{-,\gamma}(u))}_{\Hminus} \leq C(\gamma^{1/2} + M^{1/2}\gamma^{1/2})\norm{\tau_D(u)}_{\Hminus}.
    \end{align*}
        \item[2.]  If $\Omega$ has at least one convex vertex, then for each $0<a<1$ there exists a constant $C_a$ such that for all $0<M <C_a^{-1}\gamma^{-a}$, and all $u\in H^1_\gamma(\Omega)$ such that  $\tau_D(u)$ is in $\XM$,
    \begin{align*}
        \norm{(S_0(\tau_{-,0}(u))-S_\gamma(\tau_{-,\gamma}(u))}_{\Hminus} \leq C_a(\gamma^{1/3} + M^{1/2}\gamma^{a/2})\norm{\tau_D(u)}_{\Hminus}.
    \end{align*}
    \item[3.]
           If $\Omega$ has smooth boundary, then there exists a constant $C$ such that for all $0<M<C^{-1}\gamma^{-1}$, and all $u\in H^1_\gamma(\Omega)$ such that $\tau_D(u)$ is in $\XM$,
         \begin{align*}
        \norm{(S_0(\tau_{-,0}(u))-S_\gamma(\tau_{-,\gamma}(u))}_{\Hminus} \leq C(\gamma^{1/2} + M\gamma)\norm{\tau_D(u)}_{\Hminus}.
    \end{align*}
    \end{enumerate}
\end{cor}

Theorem \ref{thm:exchange} and Corollary \ref{cor:DtN} provide quantitative estimates between the local operators appearing in classical Optimized Schward Methods and their non-local counterparts, in the case of a decomposition into two sub-domains. We will see in Section \ref{sec:Claeys} that in fact our estimates can be generalized to domain decompositions with an arbitrary number of components, including those with cross-points.

\section*{Acknowledgements}

The authors thank Xavier Claeys and Euan Spence for introducing them to the problem and helpful conversations throughout the work.  The initial suggestion to work on this problem came when the first and third authors were attending the workshop ``At the Interface between Semiclassical Analysis and Numerical Analysis of Wave Scattering Problems" at the Mathematisches Forschungsinstitut Oberwolfach in October 2022.    
Y.C. was supported by NSF CAREER Grant DMS-2045494.  J.L.M. acknowledges support from NSF FRG grant DMS-2152289 and NSF Applied Math Grant DMS-2307384. 

\section{Estimates on the DtN operator} \label{sec:DtN}

In this section, we will prove the estimates on the Dirichlet-to-Neumann operator $T_\gamma$ given in Theorem \ref{thm:DtN}.  
To begin, we introduce some necessary notation that we will use throughout the proofs in this section:

 We let $\lambda_j$ be the eigenvalues of $T_{\gamma}$, with corresponding orthogonal eigenfunctions $h_j\in H^{1/2}_{\gamma}(\pa\Omega)$:
 \begin{equation}\label{eqn:hj}
 T_\gamma h_j=\lambda_j h_j, \qquad  \norm{h_j}_{L^2(\pa\Omega)} = 1.
 \end{equation}
 Note that this ensures that the functions $h_j$ are also orthogonal in $L^2(\pa\Omega)$ and $H^{-1/2}_{\gamma}(\pa\Omega)$. In particular, 
 $$
 \norm{h_j}_{H^{1/2}_{\gamma}(\pa\Omega)} = \gamma^{-1/2}\lambda_j^{1/2}, \qquad \norm{h_j}_{H^{-1/2}_{\gamma}(\pa\Omega)} = \gamma^{1/2}\lambda_j^{-1/2}.
 $$

 For each $h\in \Hplus$, we write $h = \sum_{j=1}^{\infty}a_j h_j$. Then, $\norm{h}^2_{L^2(\pa\Omega)}=\sum_{j=1}^{\infty}|a_j|^2$, and
    \begin{align*}
        \norm{h}^2_{\Hplus}  = \gamma^{-1}\sum_{j=1}^{\infty}\lambda_j|a_j|^2, \qquad    \norm{h}^2_{\Hminus}  = \gamma\sum_{j=1}^{\infty}\lambda_j^{-1}|a_j|^2.
    \end{align*}
 Given $h\in \Hplus$, we define the Rayleigh quotient by
    \begin{align} \label{eqn:Rayleigh}
        R[h] = \inf\left\{ \frac{\norm{F}^2_{\Hone}}{\gamma^{-1}\norm{h}^2_{\Lpa}} :\; \, F|_{\pa\Omega} = h\right\}.
    \end{align}
    Then, $R[h_j] = \lambda_j$, and for $h = \sum_{j=1}^{\infty}a_j h_j$,
    \begin{align*}
        R[h]% = \frac{\sum_{j=1}^{\infty}\lambda_j|a_j|^2}{\norm{h}^2_{\Lpa}}
        = \gamma \frac{\norm{h}^2_{\Hplus}}{\norm{h}^2_{\Lpa}}.
    \end{align*}

With these spectral and variational conventions defined, we  establish a series of estimates, culminating in the proof of Theorem \ref{thm:DtN}. The first is an estimate on the Steklov eigenvalues of $\Omega$, where we split into two cases depending on whether $\Omega$ has at least one convex vertex. To do this, we use estimates on Robin eigenvalues with large Robin parameter from \cite{Kh} and \cite{pank2013}.

\begin{prop} \label{prop:Steklov-poly}
Let $\Omega$ be a smooth curvilinear polygon in $\R^2$. There exist constants $c,C>0$ such that the following holds.
\begin{itemize}
\item If $\Omega$ has no convex vertices and  $0<\gamma<c$, then $\lambda_j \geq 1-C\gamma$ for $j\geq 1$.
\item If $\Omega$ has at least one convex vertex, then there exists $N\in\mathbb{Z}_{>0}$ such that for all  $0<\gamma<c$, $$\lambda_j \geq 1 - C\gamma^{2/3}, \qquad j\geq N+1.$$

In addition,  there exist $\lambda_1^*, \dots, \lambda_N^* \in (0,1)$ such that 
\begin{align*}
  |\lambda_j - \lambda_j^*| \leq C\gamma^{2/3},  \qquad 1 \leq j \leq N.
\end{align*}
If $\Omega$ has $K$ convex vertices with interior angles $\{2\alpha_k\}_{k=1}^K$,
 the collection $\{\lambda_j^*\}_{j=1}^N$ is the union of the discrete Steklov eigenvalues, with respect to the operator $-\Delta+1$, of the infinite sectors $U(\alpha_k)$ with interior angles $2\alpha_k$.  
 
 Finally, let $\lambda_j^*$ be a Steklov eigenvalue associated to the sector  $U(\alpha_{k_j})$ and let $v_j$, $V_j$ satisfy $(-\Delta+\gamma^{-2})V_j=0$ in $U(\alpha_{k_j})$, $V_j = v_j$ on $\pa U(\alpha_{k_j})$, $\gamma\pa_\nu V_j = \lambda^*_j$ on $\pa U(\alpha_{k_j})$, with $v_j$ $L^2(\pa U(\alpha_{k_j}))$-normalized. Assuming that the vertex of the sector is placed at $0$, then
\begin{align} \label{eqn:exp-decay}
    \int_{\pa U(\alpha_{k_j})}|v_j(y)|^2e^{c|y|/\gamma} d\sigma(y)+  \int_{U(\alpha_{k_j})}(|V_j(x)|^2+|\nabla V_j(x)|^2)e^{c|x|/\gamma}dx< C.
\end{align}
\end{itemize}
\end{prop}

\begin{proof1}{Proposition \ref{prop:Steklov-poly}}
For fixed $w>0$, we let $\mu_k = \mu_k(w)$ be the Robin eigenvalues of $-\Delta$, with Robin parameter $w>0$. That is, there exists a function $\phi_k = \phi_k(w)\in H^1(\Omega)$ such that
\begin{align*}
    (-\Delta - \mu_k)\phi_k =0 \text{ in }\Omega, \quad \pa_{\nu}\phi_k = w\phi_k \text{ on }\pa\Omega.
\end{align*}
There exists a constant $C>0$ such that these Robin eigenvalues have the following properties: Using Sobolev trace estimates (see Theorem 1.5.1.10 in \cite{grisvard}),
\begin{align} \label{eqn:Robin1}
    \mu_k(w)\geq -Cw^2
\end{align}
for all $k\geq1$, and $w$ sufficiently large. Moreover, from Theorem 1.2 in \cite{Kh}, there exists $N\geq0$, such that, for $1\leq k \leq N$, and $w$ sufficiently large,
\begin{align} \label{eqn:Robin2}
|\mu_k(w) - \mu_k^* w^2| \leq Cw^{4/3},
\end{align}
where $\mu_k^*<-1$ are the Robin eigenvalues with $w=1$ of the sectors $U(\alpha_m)$ with interior angle $2\alpha_m$ equal to that of a convex vertex of $\Omega$. Note that $N=0$ if $\Omega$ is smooth, or only has non-convex vertices. Also, from Proposition 5.1 in \cite{Kh}, for $k\geq N+1$, 
\begin{align} \label{eqn:Robin3}
    \mu_k(w) \geq -w^2- Cw^{4/3}
\end{align}
for $w$ sufficiently large. If $\Omega$ is smooth or has no convex vertices, then instead using the estimates from \cite{pank2013}, gives
\begin{align} \label{eqn:Robin3-smooth}
    \mu_k(w) \geq -w^2- Cw
\end{align}
for $w$ sufficiently large. 

As a consequence of \eqref{eqn:Robin1}, $\mu_k(w)$ being strictly monotonically decreasing in $w$ for fixed $k$, and since $\lim_{w\to\infty}\mu_k(w)=-\infty$, there exists $\gamma_0>0$ such that for each $\gamma$ with $0<\gamma<\gamma_0$ and each $k\geq1$, there exists a unique $w = w_k(\gamma)$ such that
\begin{align*}
    \mu_k(w) = -\gamma^{-2}.
\end{align*}
Moreover, as $\mu_k(w) \leq \mu_{k+1}(w)$, the sequence $w_k(\gamma)$ is increasing in $k$. Therefore, setting $\tau_k = \tau_k(\gamma) := \gamma w_k(\gamma)$, we have
\begin{align*}
     (-\Delta +\gamma^{-2})\phi_k =0 \text{ in }\Omega, \quad \pa_{\nu}\phi_k = \gamma^{-1}\tau_k\phi_k \text{ on }\pa\Omega,
\end{align*}
and so $\tau_k$ is an eigenvalue of the DtN operator $T_{\gamma}$, with corresponding eigenfunction $\phi_k|_{\pa\Omega}\in H^{1/2}_{\gamma}(\pa\Omega)$. Conversely, if $\lambda_j = \lambda_j(\gamma)$ is an eigenvalue of $T_{\gamma}$ with eigenfunction $h_j\in H^{1/2}(\pa\Omega)$, then defining $w_j = w_j(\gamma) = \gamma^{-1}\lambda_j(\gamma)$, we have
\begin{align*}
     (-\Delta +\gamma^{-2})H_j =0 \text{ in }\Omega, \quad \pa_{\nu}H_j = wh_j \text{ on }\pa\Omega,
\end{align*}
for a function $H_j\in H^1(\pa\Omega)$, with $H_j|_{\pa\Omega} = h_j$. Therefore, $-\gamma^{-2}$ is a Robin eigenvalue of $-\Delta$, with Robin parameter $w$, and so must be equal to $\mu_k(w)$ for some $k$. As a consequence, for $0<\gamma<\gamma_0$, the sequence  $\{\tau_k(\gamma)\}$ coincides with the sequence $\{\lambda_j(\gamma)\}$ of eigenvalues of $T_{\gamma}$.

Using $\lambda_k(\gamma) = \tau_k(\gamma) = \gamma w_k(\gamma)$, with $\gamma = (-\mu_k(w))^{-1/2}$, when $\Omega$ has convex vertices, the estimates in \eqref{eqn:Robin2} and \eqref{eqn:Robin3} give, for $\gamma>0$ sufficiently small, 
\begin{align} \label{eqn:Robin4}
    |\lambda_j(\gamma) - \lambda_j^*|\leq  C\gamma^{2/3}, 
\end{align}
for $1\leq j \leq N$, with $\lambda_j^* = (-\mu_j^*)^{-1/2}<1$, and
\begin{align} \label{eqn:Robin5}
    \lambda_j(\gamma) \geq 1 - C\gamma^{2/3},
\end{align}
for $j\geq N+1$. If instead $\Omega$ is smooth or has only non-convex vertices, the estimate in \eqref{eqn:Robin3-smooth} gives, for $\gamma>0$ sufficiently small, 
\begin{align} \label{eqn:Robin6}
    \lambda_j(\gamma) \geq 1 - C\gamma,
\end{align}
for $j\geq1$. The estimates in \eqref{eqn:Robin4}, \eqref{eqn:Robin5} and \eqref{eqn:Robin6} are the bounds on $\lambda_j$ in the statement of the proposition.

We are left to prove the exponential decay estimates on the Steklov eigenfunctions of the infinite sectors. Since the infinite sector is dilation invariant, we can rescale to $\gamma=1$, and then the required estimates follow from the decay of the Robin eigenfunctions of infinite sectors from Theorem 5.1 in \cite{khalile2017}.
\end{proof1}

For functions $h\in \XM$, we will use the following estimate on how concentrated $h$ can be near a vertex of $\Omega$. In what follows, $D_{r}(v)$ denotes the disc of radius $r$ centered at  $v$.

\begin{lem} \label{lem:vertex}
Let $\Omega$ be a smooth curvilinear polygon in $\R^2$. There exist constants $c,C>0$ and for each $0<a<1$, $0<\gamma<c$, $M>0$, $h\in \XM$, and every vertex $v$ of $\Omega$,
\begin{align*}
  \norm{h}_{L^{\infty}(\pa\Omega)} \leq C(1+M^{1/2})\norm{h}_{L^2(\pa\Omega)}, \qquad  \int_{\pa\Omega\cap D_{\gamma^{a}}(v)}h^2 \leq \min\{ C(1+M)\gamma^{a},1\}\norm{h}_{L^2(\pa\Omega)}^2.
\end{align*}

\end{lem}
\begin{proof1}{Lemma \ref{lem:vertex}}  The upper bound on the integral follows immediately from the $L^{\infty}(\pa\Omega)$ estimate in the statement, and so we just need to prove the first inequality. By dividing by $\norm{h}_{L^2(\pa\Omega)}$, we assume that $\norm{h}_{L^2(\pa\Omega)}=1$. Let $L^*$ be the maximum of $|h|$, attained at $x^*\in \pa\Omega$. Letting $w(t)$, $t\in I$, be an arc-length parameterization of $\pa\Omega$, with $x^*=w(0)$, we write $g(t) = h(w(t))$. Then,
\begin{align*}
    |g(t)-g(0)| = \left|\int_{0}^{t}g'(s)\,ds\right| \leq t^{1/2}\left(\int_{0}^{t}|g'(s)|^2\,ds\right)^{1/2} \leq Ct^{1/2}M,
\end{align*}
where we have used that $h\in \XM$. Therefore, there exists $c>0$ such that
\begin{align*}
|g(t)| \geq \tfrac{1}{2}L^* \text{ for } 0\leq t \leq c \ \min\{1,M^{-2}(L^*)^2\},
\end{align*}
and so
\begin{align*}
    \tfrac{1}{4}(L^*)^2c \ \min\{1,M^{-2}(L^*)^2\} \leq \int_I g(t)^2 \,dt \leq C.
\end{align*}
Rearranging this inequality gives the desired upper bound $\norm{h}_{L^{\infty}(\pa\Omega)} =L^* \leq C(1+M^{1/2})$ for an increased constant $C$.
\end{proof1}

In the case where $\Omega$ has at least one convex vertex, we use Lemma \ref{lem:vertex} to bound the inner product between a function $h\in \XM$ and an eigenfunction $h_j$ of $T_\gamma$, with $1\leq j\leq N$, and $N$ as in the statement of Proposition \ref{prop:Steklov-poly}.
\begin{prop}\label{prop:eigenfunction-poly}
Let $\Omega$ be a smooth curvilinear polygon in $\R^2$ with at least one convex vertex.
    There exist constants $c$, $C>0$ such that the following holds. If $N$ is as in Proposition \ref{prop:Steklov-poly}, then for $0<\gamma<c$, $M>0$, $h\in \XM$, and $1\leq j \leq N$, 
    \begin{align*}
       \left| \int_{\pa\Omega}h h_j \right| \leq C(1+M^{1/2}) \gamma^{1/2}\norm{h}_{L^2(\pa\Omega)},
    \end{align*}
    where $h_j$ is as in \eqref{eqn:hj}.
\end{prop}
\begin{proof1}{Proposition \ref{prop:eigenfunction-poly}}
We will prove this proposition by first comparing the eigenfunctions $h_j$ with the Steklov eigenfunctions of infinite sectors of a convex angle, and then use the $L^{\infty}(\pa\Omega)$ estimate from Lemma \ref{lem:vertex}.

Using the notation of Proposition \ref{prop:Steklov-poly}, fix $1\leq j \leq N$ and let $U(\alpha_{k_j})$ denote the  infinite sector that has Steklov eigenvalue $\lambda_j^*$  and  let $v_j$ be the corresponding eigenfunction of $T_\gamma$ in $U(\alpha_{k_j})$. If the function $V_j$ is the extension of $v_j$ from $\pa U(\alpha_{k_j})$ to $U(\alpha_{k_j})$, satisfying $(-\Delta+\gamma^{-2})V_j=0$,  we have the Rayleigh quotient
\begin{align} \label{eqn:poly0}
    R[v_j] = \frac{\int_{U(\alpha_{k_j})}|\nabla V_j|^2 + \gamma^{-2}V_j^2}{\gamma^{-1}\int_{\pa U(\alpha_{k_j})}v_j^2} =\gamma\int_{U(\alpha_{k_j})}|\nabla V_j|^2 + \gamma^{-2}V_j^2 =  \lambda_j^*.
\end{align}
We now define a change of variables $S_j$ to translate, rotate, and straighten the sides of a neighborhood of a convex vertex of $\Omega$ with interior angle $2\alpha_{k_j}$ to coincide with part of $U(\alpha_{k_j})$: We fix $\tfrac{2}{3}<a<1$, and work in a disc centered at this vertex of radius $\gamma^{a}$. We can rotate and translate so that, without loss of generality, the vertex is at the origin, $(s,t)=(0,0)$, and this part of the boundary of $\pa \Omega$ can be written as $t = f_{+}(s)$ for $s>0$ and $t = f_{-}(s)$ for $s<0$, where $f_\pm$ are two positive functions. Moreover, the rotation can be chosen to ensure that $\lim_{s\to0^{\pm}}f_{\pm}(s)/s = a_{\pm}$, with
\begin{align}  \label{eqn:Jacobian1}
    c \leq |a_{\pm}| \leq C, \quad  |f_{\pm}(s)/s - a_{\pm}| \leq C|s|, \quad |f'_{\pm}(s)-a_{\pm}| \leq C|s|,
\end{align}
for constants $c$, $C>0$. Denoting $U(\alpha_{k_j},\gamma^a)$ to be the part of the infinite sector $U(\alpha_{k_j})$ within a distance of $\gamma^a$ from the vertex, we can then define the function
\begin{align*}
    S_j:  U(\alpha_{k_j},\gamma^a)  \to \Omega, \qquad 
    S_j(s,t)
    =\begin{cases}
    \left(s,\frac{f_{+}(s)}{a_{+}s}t\right) & s>0,\\
    \left(s,\frac{f_{-}(s)}{a_{-}s}t\right) &s<0.
    \end{cases}
\end{align*}
which is a Lipschitz homeomorphism onto its image (and a diffeomorphism away from $s=0$). The Jacobian of this change of variables is given by the determinant of 
\begin{align*}
   J(s,t)=  \begin{pmatrix}
    1 & t\displaystyle{\frac{sf'_{\pm}(s)-f_{\pm}(s)}{a_{\pm} s^2}}\\
   0 & \displaystyle{\frac{f_{\pm}(s)}{a_{\pm}s}}
    \end{pmatrix}.
\end{align*}
In particular, since we have $|s|,|t|\leq \gamma^a$, by the bounds in \eqref{eqn:Jacobian1} we have
\begin{align} \label{eqn:Jacobian2}
|J(s,t)-Id| \leq C\gamma^a.
\end{align}
We now define $$w_j(x) = \chi_j(x)v_j(S_j^{-1}(x)),$$ with $\chi_j(x)$ a cut-off function to a neighborhood of size $\gamma^a$ about the convex vertex of angle $2\alpha_{k_j}$ of $\pa\Omega$.  Moreover, $\chi_j(x)$ is chosen so that $\norm{w_j}_{L^2(\pa\Omega)} =1.$ The function $w_j$ can be extended to $\Omega$ by $W_j(x) = \chi_j(x)V_j(S_j^{-1}(x))$, so that
\begin{align*}
    R[w_j] \leq \frac{\int_{\Omega}|\nabla W_j|^2 + \gamma^{-2}W_j^2}{\gamma^{-1}\int_{\pa \Omega}w_j^2}.
\end{align*}
Using the exponential decay properties of $v_j$ and $V_j$ from \eqref{eqn:exp-decay} in Proposition \ref{prop:Steklov-poly}, \eqref{eqn:poly0}, and the Jacobian bounds from \eqref{eqn:Jacobian2}, thus ensures that
\begin{align}\label{eqn:ry}
    R[w_j] \leq \lambda_j^* + C\gamma^a.
\end{align}
To see this we use \eqref{eqn:poly0} and that \eqref{eqn:exp-decay} and \eqref{eqn:Jacobian2} implies that
\begin{align} \label{eqn:exp-decay1}
    \left|1-\int_{\Omega\cap\text{supp}(\chi_j)}v_j(S_j^{-1}(x))^2  \right| \leq C\gamma^a,\quad     \left|1-\int_{\Omega\cap\{\chi_j(x) = \max \chi_j\}}v_j(S_j^{-1}(x))^2  \right| \leq C\gamma^a,
\end{align}
and the analogous estimate for $V_j$ and $\nabla V_j$.

The upper bound from \eqref{eqn:ry} together with the estimate on $|\lambda_j-\lambda_j^*|$ from Proposition \ref{prop:Steklov-poly}, since $\tfrac{2}{3}<a<1$, implies that
\begin{align} \label{eqn:poly1}
    R[w_j]-\lambda_j \leq C\gamma^{2/3}. 
\end{align}    
Moreover, using again the exponential decay properties from \eqref{eqn:exp-decay}, the orthogonality of the $v_j$'s, and the fact that on the support of $w_j$ we can write $w_j(x) = v_j(S_j^{-1}(x)) + (1-\chi_j(x))v_j(S_j^{-1}(x))$, gives
\begin{align} \label{eqn:poly1a}
    \qquad \qquad  \left|\int_{\pa\Omega} w_jw_k\right| \leq C\gamma^a
\end{align}
for $1\leq k \leq N$, $ k\neq j$. 

We write $w_j = \sum_{\ell=1}^{\infty} a_{j,\ell }h_\ell$, so that $\norm{w_j}^2_{L^2(\pa\Omega)} = \sum_{\ell=1}^{\infty}|a_{j,\ell}|^2=1$ and $R[w_j] = \sum_{\ell=1}^{\infty}\lambda_\ell|a_{j,\ell}|^2$.  We will show that
\begin{align} \label{eqn:poly2}
A_j+B_j:= \sum_{\ell<j:\lambda^*_{\ell}\neq \lambda^*_j}|a_{j,\ell}|^2+\sum_{\ell>j:\lambda^*_{\ell}\neq \lambda^*_j}|a_{j,\ell}|^2 \leq C\gamma^{2/3}.
\end{align}
For $j=1$, we have $A_1=0$, and using \eqref{eqn:poly1} we have
\begin{align*}
\lambda_1(1-B_1) + \lambda_{k_1}B_1 \leq R[w_1] = \sum_{\ell=1}^{\infty}\lambda_\ell|a_{1,\ell}|^2  \leq \lambda_1 + C\gamma^{2/3},
\end{align*}
where $\ell = k_1$ is the smallest index for which $\lambda_\ell^* > \lambda_1^*$. Rearranging this inequality implies \eqref{eqn:poly2} for $j=1$. The same argument also works for any $j<k_1$. In particular, using the almost orthogonality inequalities from \eqref{eqn:poly1a}, and since $\tfrac{2}{3}<a<1$, for $1\leq j<k_1$, we must have coefficients $b_{\ell,j}$ with $1-C\gamma^{2/3}\leq\sum_{\ell}|b_{\ell,j}|^2 \leq 1+C\gamma^{2/3}$, such that
\begin{align} \label{eqn:poly3}
    \norm{h_j-\sum_{1\leq \ell \leq N:\lambda_\ell^* = \lambda_j^*}b_{\ell,j}w_\ell}_{L^2(\pa\Omega)} \leq C\gamma^{2/3}.
\end{align}
For $j\geq k_1$, we first use \eqref{eqn:poly3} and the almost orthogonality inequalities from \eqref{eqn:poly1a} to obtain $A_j \leq C\gamma^{2/3}$. Moreover, using the bound on the Rayleigh quotient from \eqref{eqn:poly1} gives
\begin{align*}
    \lambda_j(1-B_j) + \lambda_{k_j}B_j - C\gamma^{2/3} \leq R[w_j] \leq \lambda_j + C\gamma^{2/3},
\end{align*}
which can be rearranged to give \eqref{eqn:poly2}. Therefore, \eqref{eqn:poly2} holds for all $1\leq j \leq N$, and hence so does \eqref{eqn:poly3}.

We can use \eqref{eqn:poly3}, together with the $L^{\infty}(\pa\Omega)$ estimate on $h$ from Lemma \ref{lem:vertex}, to complete the proof of the proposition. We normalize $h$ so that $\norm{h}_{L^2(\pa\Omega)} = 1$. Then,
\begin{align*}
    &  \left| \int_{\pa\Omega}h h_j \right|  \leq \left| \int_{\pa\Omega}h \sum_{1\leq \ell \leq N:\lambda_\ell^* = \lambda_j^*}b_{\ell,j}w_\ell \right| + \left| \int_{\pa\Omega}h \left(h_j-\sum_{1\leq \ell \leq N:\lambda_\ell^* = \lambda_j^*}b_{\ell,j}w_\ell\right) \right| \\
     &\quad\leq C(1+M^{1/2})\int_{\pa\Omega}\left|\sum_{1\leq \ell \leq N:\lambda_\ell^* = \lambda_j^*}b_{\ell,j}w_\ell\right| + C\gamma^{2/3} \\
     &\quad\leq C(1+M^{1/2})\left(\int_{\pa\Omega}e^{c\min_v\text{dist}(x,v)/\gamma}\sum_{1\leq \ell \leq N:\lambda_\ell^* = \lambda_j^*}|b_{\ell,j}|^2|w_\ell|^2\right)^{1/2}\left(\int_{\pa\Omega}e^{-c\min_v\text{dist}(x,v)/\gamma}\right)^{1/2} + C\gamma^{2/3} \\
     &\quad\leq C(1+M^{1/2})\gamma^{1/2},
\end{align*}
where to obtain the last inequality, we first used that, from \eqref{eqn:poly3}, the $b_{\ell,j}$ coefficients are bounded by a constant. Recalling that $w_\ell(x) = \chi_\ell(x)v_\ell(S_\ell^{-1}(x)),$ we then used \eqref{eqn:exp-decay} from Proposition \ref{prop:Steklov-poly} to bound the integral. This completes the proof of the proposition.
\end{proof1}

Lemma \ref{lem:vertex} also allows us to bound the Rayleigh quotient from \eqref{eqn:Rayleigh} for functions $h\in \XM$.
\begin{prop} \label{prop:Rayleigh-poly}
Let $\Omega$ be a smooth curvilinear polygon in $\R^2$, and let $h\in \XM$, with $0<M\leq \gamma^{-1}$.
If $\Omega$ has  no convex vertices,
    there exist constants $c$, $C>0$ such that for all $0<\gamma<c$,
    \begin{align*}
        R[h] \leq 1+ C(1+ M)\gamma .
    \end{align*}
    If $\Omega$ has at least one convex vertex, then, for each $a\in (0,1)$, there exist  constants $c$, $C_a>0$ such that, for all $0<\gamma<c$,
\begin{align*}
    R[h] \leq 1 + C_aM\gamma^a + C\gamma.
\end{align*}
Finally, if $\Omega$ has smooth boundary, there exist constants $c$, $C>0$ such that for all $0<\gamma<c$,
\begin{align*}
    R[h] \leq 1 + C\gamma + CM^2\gamma^2.
\end{align*}
\end{prop}
\begin{proof1}{Proposition \ref{prop:Rayleigh-poly}}  Without loss of generality, we assume that $\norm{h}_{L^2(\pa\Omega)}=1$.

Given $\gamma>0$ and $0<a<1$, let $U_{\gamma^a}$ be the subset of $\pa\Omega$, with a $\gamma^a$ neighborhood of each convex vertex of $\pa\Omega$ removed, and a $\gamma$ neighborhood of each non-convex vertex removed. We then let $w_\ell(s)$ be an arc-length parameterization of a connected component $\ell$ of $U_{\gamma^a}$ for $s\in I_{\ell,\gamma,a}$. Then, we define
\begin{align*}
    F_{\ell}:I_{\ell,\gamma,a}\times (0,\tau^*\gamma^a] \to \Omega, \qquad 
   F_{\ell} (s,t)=  w_{\ell}(s) + t\textbf{n},
\end{align*}
where, for ease of notation, we have suppressed the dependence of $F_{\ell}$ on $\gamma$ and $a$. Here, $\textbf{n}$ is the unit, inward pointing, normal to $\pa\Omega$ at a point on $U_{\gamma^a}$,  and the cut-distance constant $\tau^*>0$, depending only on the geometry of $\pa\Omega$, ensures that $F_{\ell}$  is a diffeomorphism onto its image in $\Omega$. In particular, for $\gamma>0$ sufficiently small, the image of $F_{\ell}$ is contained in $\Omega$ even when $\ell$ is part of a side of $\pa\Omega$ connected to a non-convex vertex.  The Jacobian, $J_{\ell}(s,t)$, of this change of variables satisfies
\begin{align} \label{eqn:Jacobian}
    |J_{\ell}(s,t) - Id| \leq Ct, 
\end{align}
for a constant $C$ depending on the curvature of the sides of $\pa\Omega$.

Next, we define smooth cut-off functions $\chi_0(t)$ and $\chi_{\gamma,a}(s)$ with the following properties. The function $\chi_0(t)$ is equal to $1$ for $t\leq \tfrac{1}{2}$, equal to $0$ for $t\geq1$, and with bounded derivative. The function $\chi_{\gamma,a}(s)$ is equal to $1$ on $I_{\ell,2\gamma,a}$, equal to $0$ outside $I_{\ell,\gamma,a}$, and its derivative is bounded  by $C\gamma^{-a}$ in an $O(\gamma^{a})$ neighborhood of the convex vertices of $\Omega$, and bounded by $C\gamma^{-1}$ in an $O(\gamma)$ neighborhood of the non-convex vertices, and zero otherwise. Writing $g_{1,\ell}(s) = h(w_{\ell}(s))$, we then extend $g_{1,\ell}(s)\chi_{\gamma,a}(s)$ to the image of $F_\ell$ in $\Omega$ by
\begin{align} \label{eqn:extension1}
    H_{1,\ell}(x) = g_{1,\ell}(s)\chi_{\gamma,a}(s)e^{-t/\gamma}\chi_0(t/(\tau^*\gamma^a)), \qquad x = w_{\ell}(s) + t\textbf{n}.
\end{align}
 Note that when $x\in \ell$, $H_{1,\ell}(x)=h(x)\chi_{\gamma,a}(s)$ as $t=0$.

The above functions $H_{1,\ell}(x)$ provide an extension of the part of the function $h(x)$ away from the vertices of $\pa\Omega$. Given a vertex $v$ of $\pa\Omega$, after a rotation, we can write the part of $\Omega$ in the disc $D_{2\gamma^a}(v)$ centered at $v$ of radius $2\gamma^a$ as $x_2> f_v(x_1)$, with $f_v(x_1)$ having bounded slope. Then, for $x=(x_1,x_2)\in D_{2\gamma^a}(v)$, for each vertex $v$, with $g_{2,v}(x_1) = h(x_1,f_v(x_1))$ we define
\begin{align} \label{eqn:extension2}
    H_{2,v}(x) = g_{2,v}(x_1)e^{-(x_2-f_v(x_1))/\gamma}\chi_0(x_2/\gamma^a)\tilde{\chi}_{\gamma,a,v}(x_1),
\end{align}
with $\tilde{\chi}_{\gamma,a,v}(x_1)$ a cut-off function to a neighborhood of $v$ to be defined below. Note that when $x=(x_1,x_2)\in \partial \Omega\cap D_{2\gamma^a}(v)$, $H_{2,v}(x)=h(x)\tilde\chi_{\gamma,a,v}(x_1)$ as $x_2=f_v(x_1)$. Therefore, we can define the cut-off functions $\tilde{\chi}_{\gamma,a,v}$ so that
\begin{align*}
  H(x):=  \sum_{\ell}H_{1,\ell}(x) + \sum_{v}H_{2,v}(x) = h(x) \text{ for all } x\in\pa\Omega,
\end{align*}
which in particular means that the derivative of $\tilde{\chi}_{\gamma,a,v}$ is bounded  by $C\gamma^{-a}$ in an $O(\gamma^{a})$ neighborhood of the convex vertices of $\Omega$, and bounded by $C\gamma^{-1}$ in an $O(\gamma)$ neighborhood of the non-convex vertices, and zero otherwise.

We now estimate the Rayleigh quotient
\begin{align*}
    R[h]\leq \frac{\int_{\Omega}|\nabla H|^2 + \gamma^{-2}H^2}{\gamma^{-1}\int_{\pa\Omega}h^2} = \gamma  \int_{\Omega} (|\nabla H|^2 + \gamma^{-2}H^2) .
\end{align*}
We claim that for each side $\ell$,
\begin{equation}\label{eqn:R1}
\gamma  \int_{\Omega}(|\nabla H_{1,\ell}|^2 + \gamma^{-2}H_{1,\ell}^2 )\leq 1 + C(1+M)\gamma + CM^2\gamma^2.
\end{equation}
for a constant $C$. First, we have
\begin{align*}
   \left| |\pa_tH_{1,\ell}(x)|^2 -\gamma^{-2}H_{1,\ell}(x)^2\right| &\leq Cg_{1,\ell}(s)^2e^{-c\gamma^{a-1}}1_{(0,\gamma^a)}(t), \\
   |\pa_sH_{1,\ell}(x)|^2 &\leq C(|g_{1,\ell}'(s)|^2+g_{1,\ell}(s)^2|\chi'_{\gamma,a}(s)|^2)e^{-2t/\gamma},
\end{align*}
for constants $c$, $C$. Using the Jacobian bounds from \eqref{eqn:Jacobian}, and $h\in \XM$, we therefore have
\begin{align} \nonumber
\bigg|    \gamma  \int_{\Omega}(|\nabla H_{1,\ell}|^2 + \gamma^{-2}H_{1,\ell}^2) -2\gamma^{-1}\int_{\Omega}H_{1,\ell}^2\bigg|
&\leq C\gamma\int_{I_{\ell,\gamma^{a}}}\int_{0}^{\infty}g_{1,\ell}(s)^2e^{-c\gamma^{a-1}}1_{(0,\gamma^a)}(t)\,dt\,ds 
 \\& \;\;\;\;\nonumber + C\gamma\int_{I_{\ell,\gamma^{a}}}\int_{0}^{\infty}g_{1,\ell}(s)^2|\chi_{\gamma,a}'(s)|^2e^{-2t/\gamma}\,dt\,ds 
  \\& \;\;\;\;\nonumber+ C\gamma\int_{I_{\ell,\gamma^{a}}}\int_{0}^{\infty}|g_{1,\ell}'(s)|^2e^{-2t/\gamma}\,dt\,ds \\
& \leq Ce^{-c\gamma^{a-1}}+ C\gamma^2\int_{I_{\ell,\gamma^{a}}}g_{1,\ell}(s)^2|\chi_{\gamma,a}'(s)|^2\,ds + CM^2\gamma^2. \label{eqn:H1bound1}
\end{align}
To bound the integral involving $|\chi'_{\gamma,a}(s)|^2$, we use Lemma \ref{lem:vertex}: We have $|\chi'_{\gamma,a}(s)|^2 \leq C\gamma^{-2a}$ in an $O(\gamma^{a})$ neighborhood of a convex vertex, while from Lemma \ref{lem:vertex}, the integral of $g_{1,\ell}(s)^2$ restricted to such a neighborhood is bounded by $C{(1+M)}\gamma^a$. Also, we have $|\chi'_{\gamma,a}(s)|^2 \leq C\gamma^{-2}$ in an $O(\gamma)$ neighborhood of a non-convex vertex, while from Lemma \ref{lem:vertex}, the integral of $g_{1,\ell}(s)^2$ restricted to such a neighborhood is bounded by $C{(1+M)}\gamma$.  Therefore, for $\gamma$ sufficiently small, the right hand side of \eqref{eqn:H1bound1} is bounded by
\begin{align} \label{eqn:H1bound1a}
 C{(1+M)}\gamma^{2-a} + C{(1+M)}\gamma+ CM^2\gamma^2 \leq  2C(1+M) \gamma + CM^2\gamma^2,
\end{align}
using also that $0<a<1$. 

Using the Jacobian bounds from \eqref{eqn:Jacobian} again also gives
\begin{align} \label{eqn:H1bound2}
    2\gamma^{-1}\int_{\Omega}H_{1,\ell}^2 \leq 2\gamma^{-1}\int_{I_{\ell,\gamma^{a}}}\int_{0}^{\infty}g_{1,\ell}(s)^2e^{-2t/\gamma}J_{\ell}(s,t)\,dt\,ds \leq (1+C\gamma)\norm{h}_{L^2(\pa\Omega)}^2 = 1+C\gamma.
\end{align}
Combining \eqref{eqn:H1bound1a} and \eqref{eqn:H1bound2} yields the claim in \eqref{eqn:R1}, for an increased constant $C$.

Next, we claim that for each vertex $v$,
\begin{equation}\label{eqn:R2}
\gamma  \int_{\Omega}|\nabla H_{2,v}|^2 + \gamma^{-2}H_{2,v}^2 
\leq 
\begin{cases} 
C_a{(1+M)}\gamma^{a} & \text{if $v$ is convex},\\
C{(1+M)}\gamma & \text{if $v$ is not convex}.
\end{cases}
\end{equation}

 To bound the contribution from $H_{2,v}$, note that 
\begin{align*}
\left|\nabla H_{2,v}(x)\right|^2  & \leq C\gamma^{-2}|H_{2,v}(x)|^2 \\
& + C\left(g_{2,v}'(x_1)^2|\tilde{\chi}_{\gamma,a,v}(x_1)|^2 + g_{2,v}(x_1)^2|\tilde{\chi}'_{\gamma,a,v}(x_1)|^2\right) e^{-2(x_2-f_v(x_1))/\gamma}\chi_0(x_2/\gamma^a)^2.
\end{align*}
Since $h\in \XM$, we have $\int g_{2,v}'(x_1)^2|\tilde{\chi}_{\gamma,a,v}(x_1)|^2 \leq M$, and using Lemma \ref{lem:vertex} gives
\begin{align*}
    \int g_{2,v}(x_1)^2\tilde{\chi}_{\gamma,a,v}(x_1)^2  \,dx_1 \leq C{(1+M)}\gamma^a, \quad    \int g_{2,v}(x_1)^2\tilde{\chi}'_{\gamma,a,v}(x_1)^2  \,dx_1 \leq C{(1+M)}\gamma^{-a}
    \end{align*}
    if $v$ is a convex vertex, and 
    \begin{align*}
    \int g_{2,v}(x_1)^2\tilde{\chi}_{\gamma,a,v}(x_1)^2  \,dx_1 \leq C{(1+M)}\gamma, \quad    \int g_{2,v}(x_1)^2\tilde{\chi}'_{\gamma,a,v}(x_1)^2  \,dx_1 \leq C{(1+M)}\gamma^{-1}
    \end{align*}
    if it is not. The claim in \eqref{eqn:R2} follows from combining these estimates.

    Last, since the support of $H_{1,\ell}(x)H_{2,v}(x)$ is contained in the same neighborhoods of the vertices as the support of $\tilde{\chi}_{\gamma,a,v}(x_1)$, we obtain the same estimates:
\begin{align} \label{eqn:R3}
     \gamma  \int_{\Omega}|\nabla H_{1,\ell}||\nabla H_{2,v}| + \gamma^{-2}H_{1,\ell}H_{2,v}
     \leq 
\begin{cases} 
C_a{(1+M)}\gamma^{a} & \text{if $v$ is convex},\\
C{(1+M)}\gamma & \text{if $v$ is not convex}.
\end{cases}
\end{align}

The estimate in the statement of the proposition follows from combining \eqref{eqn:R1}, \eqref{eqn:R2}, and \eqref{eqn:R3}. Here, in the case that $\Omega$ has at least one convex vertex, we are using that 
$M^2\gamma^2 < M\gamma^a$ for $0<a<1$, $M\leq\gamma^{-1}$. If $\Omega$ has no convex vertices, then we may fix a specific value of $a$,  and so we get the improved estimate, {also using $M^2\gamma^2 \leq M\gamma$ for $M\leq \gamma^{-1}$.}

Finally, if $\Omega$ has a smooth boundary, then, for a cut-distance constant $\tau^*>0$, depending on $\pa\Omega$, we can instead extend $h$ by the function $H(x) = h(w(s))e^{-t/\gamma}\chi_0(t/\tau^*)$. Here $w(s)$ an arc-length parameterization of the whole of $\pa\Omega$. Then, following the proof of \eqref{eqn:H1bound1a} and \eqref{eqn:H1bound2} above, {because it does not require the cut-off function $\chi_{\gamma,a}$}, gives the upper bound of
    \begin{align*}
        R[h] \leq 1 + C\gamma + CM^2\gamma^2,
    \end{align*}
    as required.
    \end{proof1}

We next use this upper bound on the Rayleigh quotient to control the coefficients of a function $h\in \XM$ when written in terms of the basis of Steklov eigenfunctions $h_j$.
\begin{prop} \label{prop:coeff-poly}
Let $\Omega$ be a smooth curvilinear polygon in $\R^2$, and let $h=\sum_{j=1}^{\infty}a_jh_j\in \XM$ with $\norm{h}_{\Lpa}^2=1$ and $0<M\leq \gamma^{-1}$.
 Then, if $\Omega$ has no convex vertices, there exist constants $c$, $C$ such that for all $0<\gamma<c$ and each $k\geq0$, 
\begin{align*}
    \sum_{j:\lambda_j\geq 1+2^{-k}}|a_j|^2 + \sum_{j:\lambda_j\geq 1+2^{-k}}\lambda_j|a_j|^2 \leq C 2^{k}(1+M)\gamma.
\end{align*}
If $\Omega$ has at least one convex vertex, then for each $0<a<1$, there exist constants $c$, $C_a$ such that for all $0<\gamma<c$ and each $k\geq0$, 
\begin{align*}
    \sum_{j:\lambda_j\geq 1+2^{-k}}|a_j|^2 + \sum_{j:\lambda_j\geq 1+2^{-k}}\lambda_j|a_j|^2 \leq C_a2^{k}(\gamma^{2/3}+M\gamma^{a}).
\end{align*}
{Finally, if $\Omega$ has smooth boundary, there exist constants $c$, $C$ such that for all $0<\gamma<c$ and each $k\geq0$, 
\begin{align*}
    \sum_{j:\lambda_j\geq 1+2^{-k}}|a_j|^2 + \sum_{j:\lambda_j\geq 1+2^{-k}}\lambda_j|a_j|^2 \leq C 2^{k}(\gamma + M^2\gamma^2).
\end{align*}}
\end{prop}
\begin{proof1}{Proposition \ref{prop:coeff-poly}}
 Define $c_k = \sum_{j:\lambda_j\geq 1+2^{-k}}|a_j|^2$. We first write
\begin{align*}
    R[h] = \sum_{j=1}^{\infty}\lambda_j|a_j|^2 & \geq \sum_{j> N: \lambda_j<1+2^{-k}}\lambda_j|a_j|^2 + \sum_{j: \lambda_j\geq 1+2^{-k}}\lambda_j|a_j|^2.
\end{align*}
When $\Omega$ has no convex vertices (including the case of a smooth boundary), we can use the lower bound on $\lambda_j$ from Proposition \ref{prop:Steklov-poly} for $j\geq N+1$, to get the lower bound
\begin{align}\label{eqn:coeff-poly1} \nonumber
    R[h] \geq (1-C\gamma)(1-c_k)+\sum_{j: \lambda_j\geq 1+2^{-k}}\lambda_j|a_j|^2 &\geq (1-C\gamma)(1-c_k)+(1+2^{-k})c_k \\ 
    &\geq 1-C\gamma + 2^{-k}c_k.
\end{align}
If instead $\Omega$ has at least one convex vertex, we use the lower bound on $\lambda_j$ from Proposition \ref{prop:Steklov-poly} and the bound on $|a_j|$ for $j\leq N$ from Proposition \ref{prop:eigenfunction-poly} to get
\begin{align} \label{eqn:coeff-poly2} \nonumber
    R[h] \geq (1-C\gamma^{2/3})(1-c_k-M\gamma) +  \sum_{j: \lambda_j\geq 1+2^{-k}}\lambda_j|a_j|^2&\geq (1-C\gamma^{2/3})(1-c_k-M\gamma) +  (1+2^{-k})c_k \\
    &\geq 1-C\gamma^{2/3}-M\gamma +2^{-k}c_k.
\end{align}
In each case, combining the estimates in \eqref{eqn:coeff-poly1} and \eqref{eqn:coeff-poly2} with the upper bound on $R[h]$ from Proposition \ref{prop:Rayleigh-poly}, implies the required upper bounds on $c_k$.

From \eqref{eqn:coeff-poly1} and \eqref{eqn:coeff-poly2}, we have
\begin{align*}
     \sum_{j: \lambda_j\geq 1+2^{-k}}\lambda_j|a_j|^2 \leq R[h]-1 + c_k + C\gamma
\end{align*}
when $\Omega$ has no convex vertices, and
\begin{align*}
    \sum_{j: \lambda_j\geq 1+2^{-k}}\lambda_j|a_j|^2 \leq R[h]-1 + C\gamma^{2/3} + CM\gamma+c_k
\end{align*}
when it does. Using the upper bound on $c_k$, and the estimates on $R[h]$ from Proposition \ref{prop:Rayleigh-poly} then imply the remaining estimates in the proposition.
\end{proof1}
Before completing the proof of Theorem \ref{thm:DtN}, we finally need an estimate relating the $H_\gamma^{-1/2}(\pa\Omega)$ and $L^2(\pa\Omega)$ norms of functions in $\XM$.
\begin{lem} \label{lem:H-1/2}
    If $\Omega$ has no convex vertices, there exist constants $c$, $C>0$ such that for all $0<\gamma<c$, $0<M<C^{-1}\gamma^{-1}$, and $h\in \XM$,
     \begin{align*}
        \norm{h}_{H^{-1/2}_{\gamma}(\pa\Omega)} \geq \tfrac{1}{2}\gamma^{1/2}\norm{h}_{L^2(\pa\Omega)}.
    \end{align*}
    If  $\Omega$ has at least one convex vertex, then for each $0<a<1$, there exist constants $c$, $C_a>0$ such that for all $0<\gamma<c$, $0<M<C_a^{-1}\gamma^{-a}$, and $h\in \XM$, this same lower bound on $ \norm{h}_{H^{-1/2}_{\gamma}(\pa\Omega)}$ holds. 
\end{lem}
\begin{proof1}{Lemma \ref{lem:H-1/2}}
For $h=\sum_{j=1}^{\infty}a_jh_j\in \XM$ with $\norm{h}_{L^2(\pa\Omega)} = 1$, we can write
 \begin{align*}
      \norm{h}^2_{\Hminus}  = \gamma\sum_{j=1}^{\infty}\lambda_j^{-1}|a_j|^2 \geq \gamma\sum_{j:\lambda_j\leq 2}\lambda_j^{-1}|a_j|^2 \geq \tfrac{1}{2}\gamma\sum_{j:\lambda_j\leq 2}|a_j|^2.
 \end{align*}
 By Proposition \ref{prop:coeff-poly}, for $\gamma$ and $M\gamma$ sufficiently small (and $M\gamma^a$ also sufficiently small in the case where $\Omega$ has a convex vertex), we have $\tfrac{1}{2}\gamma\sum_{j:\lambda_j\leq 2}|a_j|^2\geq \tfrac{1}{4}\gamma$, and this completes the proof of the lemma.
\end{proof1}

We now prove Theorem \ref{thm:DtN}.
\begin{proof1}{Theorem \ref{thm:DtN}}
Let $h=\sum_{j=1}^{\infty}a_jh_j\in \XM$, with $\norm{h}_{\Lpa}=1$. Then,
\begin{align*}
 \norm{(T_{\gamma}-Id)h}_{\Hminus}^2 =   \norm{\textstyle\sum_{j=1}^\infty (\lambda_j-1)a_j}_{\Hminus}^2 = \gamma\sum_{j=1}^{\infty} \lambda_j^{-1}(\lambda_j-1)^2|a_j|^2.
\end{align*}
We decompose this sum as
\begin{align}
     \norm{(T_{\gamma}-Id)h}_{\Hminus}^2=
    &\,\gamma\sum_{j:\lambda_j\leq 1} \lambda_j^{-1}(\lambda_j-1)^2|a_j|^2 \notag\\&+ \gamma\sum_{k=0}^{\infty}\sum_{j:1+2^{-k-1} \leq \lambda_j < 1+2^{-k}} \lambda_j^{-1}(\lambda_j-1)^2|a_j|^2 + \gamma\sum_{j:\lambda_j\geq 2} \lambda_j^{-1}(\lambda_j-1)^2|a_j|^2. \label{eqn:DtN-thm1}
\end{align}
We now split into three cases:

1. If $\Omega$ has no convex vertices, then we use Proposition \ref{prop:Steklov-poly} to bound the first term in \eqref{eqn:DtN-thm1}
by $C\gamma^{3}$ and Proposition \ref{prop:coeff-poly} with $k=0$ to bound the third term by $C\gamma^2(1+M)$. Using Proposition \ref{prop:coeff-poly} for all $k\geq 1$, we can bound the second term in \eqref{eqn:DtN-thm1} by
\begin{align*}
    C\gamma\sum_{k=0}^{\infty} 2^{-2k}2^{k}\gamma(1+M) \leq C\gamma^2(1+M).
\end{align*}
Thus,
\begin{align*}
 \norm{(T_{\gamma}-Id)h}_{\Hminus}^2 \leq C\gamma^2(1+M),
 \end{align*}
 and to finish the proof we apply Lemma \ref{lem:H-1/2}.

2. If instead, $\Omega$ has at least one convex vertex, then we now use Propositions \ref{prop:Steklov-poly} and \ref{prop:eigenfunction-poly} to this time bound the first term in \eqref{eqn:DtN-thm1}
by $C(\gamma^{7/3}+\gamma^2M)$ and Proposition \ref{prop:coeff-poly} with $k=0$ to bound the third term by $C_a\gamma(\gamma^{2/3}+M\gamma^a)$. Again, using Proposition \ref{prop:coeff-poly} for all $k\geq 1$, bounds the second term in \eqref{eqn:DtN-thm1} by
\begin{align*}
    C_a\gamma\sum_{k=0}^{\infty} 2^{-2k}2^{k}(\gamma^{2/3}+M\gamma^{a}) \leq C_a\gamma(\gamma^{2/3}+M\gamma^{a}).
\end{align*}
Thus,
\begin{align*}
 \norm{(T_{\gamma}-Id)h}_{\Hminus}^2 \leq C_a\gamma(\gamma^{2/3}+M\gamma^{a}),
 \end{align*}
 and to finish the proof we again apply Lemma \ref{lem:H-1/2}.

3. Finally, if $\Omega$ has smooth boundary, then we again use Proposition \ref{prop:Steklov-poly} to bound the first term in \eqref{eqn:DtN-thm1}
by $C\gamma^{3}$ and Proposition \ref{prop:coeff-poly} with $k=0$ to bound the third term by $C\gamma(\gamma+M^2\gamma^2)$. Using Proposition \ref{prop:coeff-poly} for all $k\geq 1$, we can bound the second term in \eqref{eqn:DtN-thm1} by
\begin{align*}
    C\gamma\sum_{k=0}^{\infty} 2^{-2k}2^{k}(\gamma+M^2\gamma^2) \leq C\gamma(\gamma+M^2\gamma^2).
\end{align*}
Thus,
\begin{align*}
 \norm{(T_{\gamma}-Id)h}_{\Hminus}^2 \leq C\gamma(\gamma+M^2\gamma^2),
 \end{align*}
 and to finish the proof we again apply Lemma \ref{lem:H-1/2}.
\end{proof1} 

\section{Estimates on the exchange operator} \label{sec:exchange}
In this section, we will prove Theorem \ref{thm:exchange}.
We continue to write  $\Omega_1$ for a bounded planar domain, with boundary given by a smooth, curvilinear polygon, and set $\Omega_0 = \R^2\backslash \bar{\Omega}_1$. We let $\Gamma = \pa\Omega_0= \pa\Omega_1$ correspond to the shared boundary.
We recall that the exchange operator,  $\Pi_\gamma$, acts on $\phi=(\phi_0,\phi_1)\in\HN{-1/2}(\Gamma)$ and is defined by \eqref{defn:exchange}.
Theorem \ref{thm:exchange} approximates $\Pi_\gamma$ by its local counterpart $\Pi_0$ defined in \eqref{eqn:Pi0}.
To prove Theorem \ref{thm:exchange}, we first write $\Pi_\gamma-\Pi_0$ in terms of an integral operator, which involves a  Neumann-Poincar\'e operator for $-\Delta+\gamma^{-2}$. We will then prove bounds for the kernel of this operator, which will be used to establish the estimates in Theorem \ref{thm:exchange}.
\begin{lem}\label{lem:Pi}
    The exchange operator can be written as
    \begin{align*}
     \Pi\begin{pmatrix} \phi_0\\\phi_1 \end{pmatrix} = \Pi_0
    \begin{pmatrix} \phi_0\\\phi_1 \end{pmatrix} + \begin{pmatrix}   -A& -A\\ A& A \\ \end{pmatrix}
    \begin{pmatrix} \phi_0\\\phi_1 \end{pmatrix},
\end{align*}
where $A$ is the operator acting on $H_{\gamma}^{-1/2}(\pa\Omega_k)$  for $k=0,1$ defined  by
\begin{align*}
    Af(x) = 2\int_{\pa\Omega_0} \left(\pa_{n_0,x}\mathcal{G}_{\gamma} \right)(x-y)f(y)\,d\sigma(y),
\end{align*}
with $\mathcal{G}_{\gamma}(x) =K_0(|x|/\gamma)$ and $K_0$ being the modified Bessel function  of the second kind of order $0$.
\end{lem}
\begin{proof1}{Lemma \ref{lem:Pi}}
By \eqref{eqn:Psi1}, $ \PsiN{}(\phi)(x) = \PsiN{0}(\phi_0)(x) + \PsiN{1}(\phi_1)(x)$ with
\begin{align}
\PsiN{0}(\phi_0)(x):=\int_{\pa\Omega_0}\mathcal{G}_\gamma(x-y)\phi_0(y)\,d\sigma(y), \qquad  \PsiN{1}(\phi_1)(x):= \int_{\pa\Omega_1}\mathcal{G}_\gamma(x-y)\phi_1(y)\,d\sigma(y).
\end{align}
Then, since $n_0 = -n_1$ on $\pa\Omega_0 = \pa\Omega_1$, and using $int$, $ext$ to signify derivatives evaluated from within the interior, exterior of $\Omega_j$, we have 
\begin{align*}
    \pa_{n_0}^{int}\PsiN{}(\phi)  = \pa_{n_0}^{int}\PsiN{0}(\phi_0) - \pa_{n_1}^{ext}\PsiN{1}(\phi_1), \qquad
     \pa_{n_1}^{int}\PsiN{}(\phi)  = -\pa_{n_0}^{ext}\PsiN{0}(\phi_0) + \pa_{n_1}^{int}\PsiN{1}(\phi_1).
\end{align*}
We now use jump formulas for single layer potentials to write
\begin{align*}
    \pa_{n_0}^{int}\PsiN{0}(\phi_0)  = \tfrac{1}{2}\phi_0 + S^{*,0}\phi_0, \qquad 
        \pa_{n_0}^{ext}\PsiN{0}(\phi_0)  = -\tfrac{1}{2}\phi_0 + S^{*,0}\phi_0.
\end{align*}
Note that the jump formulas are analogous to those for the single layer potential for $\Delta$, as for instance described in \cite{folland2020introduction}, Chapter $3.E$, because $\mathcal{G}_{\gamma}(x)$ has the same logarithmic singularity at $x=0$ for any $\gamma>0$.
Here, the operator $S^{*,0}$ acts on $H_{\gamma}^{-1/2}(\pa\Omega_0)$ and is defined by
\begin{align*}
    S^{*,0}f(x) = \int_{\pa\Omega_0} \left(\pa_{n_0,x}\mathcal{G}_{\gamma} \right)(x-y)f(y)\,d\sigma(y).
\end{align*}
We have the same equations for the subscript $1$, and since $n_0=-n_1$, we have $S^{*,0} = -S^{*,1}$. 

Therefore, 
\begin{align*}
    \pa_{n_0}^{int}\PsiN{}(\phi) = \pa_{n_0}^{int}\PsiN{0}(\phi_0) - \pa_{n_1}^{ext}\PsiN{1}(\phi_1)   = \tfrac{1}{2}\phi_0 +S^{*,0}\phi_0 + \tfrac{1}{2}\phi_1 - S^{*,1}\phi_1,
\end{align*}
and the analogous expression for   $\pa_{n_1}^{int}\PsiN{}(\phi)$ (with the $0$ and $1$ subscripts reversed). In matrix notation, \eqref{defn:exchange} yields
\begin{align*}
    \Pi\begin{pmatrix} \phi_0\\\phi_1 \end{pmatrix} &  = \begin{pmatrix}\phi_0 -\phi_0 - 2S^{*,0}\phi_0 -\phi_1 + 2S^{*,1}\phi_1 \\ \phi_1 -\phi_1 - 2S^{*,1}\phi_1 -\phi_0 + 2S^{*,0}\phi_0  
    \end{pmatrix}\\
    & = \begin{pmatrix}0 &-1  \\-1 &0 \\ \end{pmatrix}
    \begin{pmatrix} \phi_0\\\phi_1 \end{pmatrix} + \begin{pmatrix}   -2S^{*,0}& 2S^{*,1}\\ 2S^{*,0}& -2S^{*,1} \\ \end{pmatrix}
    \begin{pmatrix} \phi_0\\\phi_1 \end{pmatrix}.
\end{align*}
Setting $A = 2S^{*,0} = -2S^{*,1}$, concludes the proof.
\end{proof1}
As shown in Corollary 5.1 of \cite{Cl}, the exchange operator satisfies various mapping properties:
\begin{lem}[Corollary 5.1 in \cite{Cl}] \label{lem:Pi-prop}
    The operator $\Pi_\gamma$ from Definition \ref{defn:exchange} has the following properties:
\begin{enumerate}
    \item[a)] $\Pi_\gamma$ is an isometry on  $\HN{-1/2}(\Gamma)$, with $\Pi_\gamma^2 = Id$.

    \item[b)] The operator $Id-\Pi_\gamma$ has kernel $\begin{pmatrix}  f\\-f\end{pmatrix}$ for any $f\in H_{\gamma}^{-1/2}(\pa\Omega_0)$. Equivalently, the kernel consists of those elements in $\HN{-1/2}(\Gamma)$ which are the Neumann traces of functions in $H^1(\R^2)$.

    \item[c)] The operator $Id+\Pi_\gamma$ has kernel $\begin{pmatrix}T_\gamma^0h \\ T_\gamma^1h \end{pmatrix}$, where $h$ is the restriction of any $H^1(\R^2)$ function to the common boundary of $\Omega_0$ and $\Omega_1$, and $T_{\gamma}^k:H^{1/2}_{\gamma}(\pa\Omega_k)\to H^{-1/2}_{\gamma}(\pa\Omega_k)$ is the Dirichlet-to-Neumann (DtN) operator of $\Omega_k$, for the operator $-\Delta+\gamma^{-2}$.
\end{enumerate}
\end{lem}
From Lemma \ref{lem:Pi} and the following corollary, we see that the exchange operator $\Pi_\gamma$ and the Dirichlet-to-Neumann operators $T_\gamma^k$ are intricately linked. 
\begin{cor} \label{cor:Pi-prop}
    As a consequence of Lemma \ref{lem:Pi-prop} c), the operator $A$ from Lemma \ref{lem:Pi} and DtN operators $T_{\gamma}^k$ satisfy
    \begin{align*}
        A(T_\gamma^0+T_\gamma^1)h = (T_\gamma^0-T_\gamma^1)h,
    \end{align*}
    for any $h$ which is the restriction of a $H^1(\R^2)$ function to the common boundary of $\Omega_0$ and $\Omega_1$.
\end{cor}
\begin{proof1}{Corollary \ref{cor:Pi-prop}}
Using the expression for $\Pi_\gamma$ from Lemma \ref{lem:Pi} gives
\begin{align*}
    Id + \Pi_\gamma = \begin{pmatrix}   1-A& -1-A\\-1+ A& 1+A \\ \end{pmatrix},
\end{align*}
and then the corollary follows from Lemma \ref{lem:Pi-prop} c).
\end{proof1}

     We note that we can use Corollary \ref{cor:Pi-prop} to write
    \begin{align*}
        Ah & = \tfrac{1}{2}A(h-T^0_\gamma h)+\tfrac{1}{2}A(h-T^1_\gamma h) + \tfrac{1}{2}A(T_\gamma^0+T_\gamma^1)h \\
        & = \tfrac{1}{2}A(h-T^0_\gamma h)+\tfrac{1}{2}A(h-T^1_\gamma) + \tfrac{1}{2}(T_\gamma^0-T_\gamma^1)h \\
        & = \tfrac{1}{2}A(h-T^0_\gamma  h)+\tfrac{1}{2}A(h-T^1_\gamma h) + \tfrac{1}{2}(T_\gamma^0h-h)+\tfrac{1}{2}(h-T_\gamma^1h).
    \end{align*}
    In particular, using the bound on $A$ coming from Lemma \ref{lem:Pi-prop} a), and Theorem \ref{thm:DtN} 2., for $h\in \XM$, we can write
    \begin{align*}
        Ah = F_0 + F_1,
    \end{align*}
    with $F_k\in H^{-1/2}_{\gamma}(\pa\Omega_k)$ satisfying the smallness estimate on the right hand side of Theorem \ref{thm:DtN} 2. In the proof of Theorem \ref{thm:exchange} below, we will obtain bounds on $Ah$ directly, using the properties of the kernel of $A$.

To prove Theorem \ref{thm:exchange}, we will use the integral operator representation of $A$ from Lemma \ref{lem:Pi} by obtaining estimates on its kernel.
\begin{prop}\label{prop:kernel-bounds}
    There exist constants $c$, $C$ such that the kernel $K_{\gamma}(x,y):=2\left(\pa_{n_0,x}\mathcal{G}_{\gamma} \right)(x-y)$ of the operator $A$ from Lemma \ref{lem:Pi} satisfies the following bounds.
    \begin{enumerate}
        \item[1.]  Let $\Omega_1$ be a smooth domain. Then, for $x\in\pa\Omega_0=\pa\Omega_1$, $y\in \bar{\Omega}_k$, $k=0,1$
        \begin{align*}
            |K_{\gamma}(x,y)| + \gamma|\nabla_yK_{\gamma}(x,y)| \leq Ce^{-c|x-y|/\gamma}.
        \end{align*}
        \item[2.] Let $\Omega_1$ be a smooth curvilinear polygon. Then, for $x,y\in\pa\Omega_0 = \pa\Omega_1$, 
         \begin{align*}
         |K_{\gamma}(x,y)|  \leq C\frac{1}{\min_{i}
         \left(\emph{dist}(x,v_i)+\emph{dist}(y,v_i)\right)}e^{-c|x-y|/\gamma},
        \end{align*}
        where $v_i$ are the vertices of $\pa\Omega_1$.
    \end{enumerate}
\end{prop}
We will first use these bounds to prove Theorem \ref{thm:exchange}, and then prove the proposition.
\begin{proof1}{Theorem \ref{thm:exchange}}

\emph{Part 1.} We first prove the theorem in the case that $\Omega_1$ is a smooth domain. We will use Proposition \ref{prop:kernel-bounds} to show that there exists a constant $C$ such that, for $k=0,1$,
\begin{align} \label{eqn:A-L2}
    \norm{Ag}_{L^2(\pa\Omega_k)} & \leq C\gamma\norm{g}_{L^2(\pa\Omega_k)}, \\ \label{eqn:A-H12}
      \norm{Af}_{H_{\gamma}^{-1/2}(\pa\Omega_k)} & \leq C\gamma\norm{f}_{H_{\gamma}^{-1/2}(\pa\Omega_k)}. 
\end{align}
Here $A$ is the operator from Lemma \ref{lem:Pi}, with $g\in L^2(\pa\Omega_k)$ and $f\in H^{-1/2}_{\gamma}(\pa\Omega_k)$. The estimates on $\Pi_\gamma-\Pi_0$ in Part 1 of Theorem \ref{thm:exchange} then follow from Lemma \ref{lem:Pi} by Sobolev interpolation.

To prove the $L^2(\pa\Omega_k)$ estimate \eqref{eqn:A-L2}, we first note that from Part 1 of Proposition \ref{prop:kernel-bounds}, there are constants $c,C>0$ such that
\begin{align*}
    \sup_{x\in\pa\Omega_k} \int_{\pa\Omega_k}|K_{\gamma}(x,y)|\,d\sigma(y)+ \sup_{y\in\pa\Omega_k} \int_{\pa\Omega_k}|K_{\gamma}(x,y)|\,d\sigma(x) \leq C\int_{-\infty}^{\infty}e^{-c|t|/\gamma}\,dt \leq C \gamma.
\end{align*}
The constant $C$ will change from line-to-line below. Using Schur's test, 
\begin{align} \nonumber
    \int_{\pa\Omega_k}|Ag(x)|^2\,d\sigma(x) & \leq \int_{\pa\Omega_k}\left(\int_{\pa\Omega_k}|K_{\gamma}(x,y)|\,d\sigma(y)\right)\left(\int_{\pa\Omega_k}|K_{\gamma}(x,y)||g(y)|^2\,d\sigma(y)\right)\,d\sigma(x)\\ \nonumber
    & \leq C\gamma\int_{\pa\Omega_k}\left(\int_{\pa\Omega_k}|K_{\gamma}(x,y)|\,d\sigma(x)\right)|g(y)|^2\,d\sigma(y)\\ \label{eqn:Schur}
& \leq C\gamma^2\int_{\pa\Omega_k}|g(y)|^2\,d\sigma(y).
\end{align}
This give the desired $A:L^2(\pa\Omega_k)\to L^2(\pa\Omega_k)$ estimate.

To prove the $A:H_{\gamma}^{-1/2}(\pa\Omega_k)\to H_{\gamma}^{-1/2}(\pa\Omega_k)$ estimate we proceed as follows. Let the operator $A^*:L^2(\pa\Omega_k)\to H^{1}_{\gamma}(\Omega_k)$ be given by
\begin{align*}
    A^*h(y) = \int_{\pa\Omega_k}K(x,y)h(x)\,d\sigma(x).
\end{align*}
Then, we will show that for $h\in H_\gamma^{-1/2}(\pa\Omega_k)$,
\begin{align} \label{eqn:dual}
    \norm{A^*h}_{H^1_{\gamma}(\Omega_k)} \leq C\gamma^{1/2}\norm{h}_{L^{2}(\pa\Omega_k)}, \qquad
    \norm{h}_{L^2(\pa\Omega_k)} \leq C\gamma^{1/2}\norm{h}_{H_{\gamma}^{1/2}(\pa\Omega_k)}.
\end{align}
Assuming these bounds, we have
\begin{align*}
    \left|\langle Af,h\rangle\right| = \left|\langle f,A^*h\rangle\right|  & \leq \norm{f}_{H_{\gamma}^{-1/2}(\pa\Omega_k)}\norm{A^*h|_{\pa\Omega_k}}_{H_{\gamma}^{1/2}(\pa\Omega_k)} \\
    & \leq \norm{f}_{H_{\gamma}^{-1/2}(\pa\Omega_k)}\norm{A^*h}_{H_{\gamma}^{1}(\Omega_k)} \\
    &\leq C\gamma^{1/2}\norm{f}_{H_{\gamma}^{-1/2}(\pa\Omega_k)}\norm{h}_{L^{2}(\pa\Omega_k)}  \leq C\gamma \norm{f}_{H_{\gamma}^{-1/2}(\pa\Omega_k)}\norm{h}_{H_{\gamma}^{1/2}(\pa\Omega_k)}, 
\end{align*}
giving the desired $A:H^{-1/2}_{\gamma}(\pa\Omega_k)\to H_{\gamma}^{-1/2}(\pa\Omega_k)$ estimate. 

We are left to prove \eqref{eqn:dual}. The second estimate in \eqref{eqn:dual} can be proved as follows: For $\Omega_k$ Lipschitz, and $v\in H^1_{\gamma}(\Omega_k)$, by Theorem 1.5.1.10 in \cite{grisvard} we have
\begin{align*}
    \int_{\pa\Omega_k}v^2\,d\sigma \leq C\gamma\left( \int_{\Omega_k}|\nabla v|^2 + \gamma^{-2}\int_{\Omega_k}v^2\right) = C\gamma \norm{v}^2_{H^1_{\gamma}(\Omega_k)}.
\end{align*}
Applying this to any $\phi_h\in H^1_{\gamma}(\Omega_k)$, with $\phi_h|_{\pa\Omega_k}=h$, gives $\norm{h}_{L^2(\pa\Omega_k)} \leq C\gamma^{1/2}\norm{\phi_h}_{H^1_{\gamma}(\Omega_k)}$. The second bound in \eqref{eqn:dual} then follows from the definition of the $H^{1/2}_{\gamma}(\pa\Omega_k)$ norm.

We are finally left to prove the first bound in \eqref{eqn:dual}. Note from Proposition \ref{prop:kernel-bounds}, Part 1, that
\begin{align*}
\sup_{x\in\pa\Omega_k}\int_{\Omega_k}|K_\gamma(x,y)|\,dy \leq \sup_{x\in\pa\Omega_k}\int_{\Omega_k}Ce^{-c|x-y|/\gamma}\,dy \leq C\gamma^2.
\end{align*}
Therefore, following Schur's test, analogously to \eqref{eqn:Schur}, gives
\begin{align} \label{eqn:dual1}
    \int_{\Omega_k}|A^*h(y)|^2\,dy \leq C\gamma^{3}\int_{\pa\Omega_k}|h(x)|^2\,d\sigma(x).
\end{align}
Using the bound on $\nabla_yK(x,y)$ from Proposition \ref{prop:kernel-bounds} 1. instead  gives
\begin{align} \label{eqn:dual2}
    \int_{\Omega_k}|\nabla A^*h(y)|^2\,dy \leq C\gamma^{3}\gamma^{-2}\int_{\pa\Omega_k}|h(x)|^2\,d\sigma(x).
\end{align}
Combining \eqref{eqn:dual1} and \eqref{eqn:dual2}, using the definition of the $H_{\gamma}^1(\Omega_k)$-norm then  gives \eqref{eqn:dual}.

\emph{Part 2.} We now prove Theorem \ref{thm:exchange} in the case that $\Omega_1$ is a smooth curvilinear polygon. We let $d(x)$ be the minimum distance from $x\in\pa\Omega_1$ to the vertices of $\Omega_1$.  We then write the kernel $K_{\gamma}(x,y)$ as
\begin{align*}
  K_{\gamma}(x,y) = \sum_{j:2^j<1/\gamma}K^j_{\gamma}(x,y) + \tilde{K}_{\gamma}(x,y),
\end{align*}
$$
K^j_{\gamma}(x,y):=1_{\{d(x)\in(2^{-j},2^{-j+1}]\}}(x)K_{\gamma}(x,y), \qquad \tilde{K}_{\gamma}(x,y)=1_{\{d(x)\leq \gamma\}}(x)K_{\gamma}(x,y),
$$
where, for example, $1_{\{d(x)\leq\gamma\}}(x)$ is the indicator function of the set $\{d(x)\leq\gamma\}$. This then breaks the operator $A$ into the sum $\sum_{j}A_j+\tilde{A}$. Using Part 2 of Proposition \ref{prop:kernel-bounds}, we have
\begin{align*}
    |K^j_{\gamma}(x,y)| \leq C2^j1_{\{d(x)\in(2^{-j},2^{-j+1}]\}}(x)e^{-c|x-y|/\gamma}
\end{align*}
for $x,y\in\pa\Omega_k$. To handle ${A}_j$, we apply the $L^{\infty}(\pa\Omega_k)$ estimate from Lemma \ref{lem:vertex}, to see that for $g\in \XM$, 
\begin{align*}
   |{A}_jg(x)| \leq  \int_{\pa\Omega_k}|{K}^j_{\gamma}(x,y)| |g(y)|\,d\sigma(y) &\leq C1_{\{d(x)\in(2^{-j},2^{-j+1}]\}}(x)\int_{\pa\Omega_k}2^je^{-c|x-y|/\gamma} |g(y)|\,d\sigma(y)\\&\leq C\gamma {(1+M^{1/2})}2^j\norm{g}_{L^2(\pa\Omega_k)}1_{\{d(x)\in(2^{-j},2^{-j+1}]\}}(x).
\end{align*}
This gives the $L^2(\pa\Omega_k)$-estimate
\begin{align} \label{eqn:exchange-poly1}
    \norm{{A}_jg}_{L^2(\pa\Omega_k)} \leq C\gamma {(1+M^{1/2})}2^{j/2}\norm{g}_{L^2(\pa\Omega_k)}.
\end{align}
 To handle $\tilde{A}$, we use the estimate from Part 2 of Proposition \ref{prop:kernel-bounds}, and again apply the $L^{\infty}(\pa\Omega_k)$ estimate from Lemma \ref{lem:vertex}, to see that for $g\in \XM$, 
\begin{align*}
   |\tilde{A}g(x)| \leq  \int_{\pa\Omega_k}|\tilde{K}_{\gamma}(x,y)| |g(y)|\,d\sigma(y) &\leq C1_{\{d(x)\leq\gamma\}}(x)\int_{\pa\Omega_k}\frac{1}{d(x) + d(y)}\gamma^{-1}e^{-c|x-y|/\gamma} |g(y)|\,d\sigma(y)\\&\leq C{(1+M^{1/2})}\norm{g}_{L^2(\pa\Omega_k)}1_{\{d(x)\leq\gamma\}}(x)\ln d(x).
\end{align*}
 This gives the $L^2(\pa\Omega_k)$-estimate
\begin{align} \label{eqn:exchange-poly2}
    \norm{\tilde{A}g}_{L^2(\pa\Omega_k)} \leq C\gamma^{1/2} {(1+M^{1/2})}\norm{g}_{L^2(\pa\Omega_k)}.
\end{align}
From \eqref{eqn:exchange-poly1} and \eqref{eqn:exchange-poly2} we have
\begin{align} \label{eqn:exchange-poly3} \nonumber
    \norm{Ag}_{L^2(\pa\Omega_k)} &\leq C\sum_{j:2^j<1/\gamma}\gamma {(1+M^{1/2})}2^{j/2}\norm{g}_{L^2(\pa\Omega_k)} + C\gamma^{1/2} {(1+M^{1/2})}\norm{g}_{L^2(\pa\Omega_k)} \\
    &\leq C\gamma^{1/2}{(1+M^{1/2})}\norm{g}_{L^2(\pa\Omega_k)}
\end{align}
for all $g\in \XM$. From Lemma \ref{lem:Pi}, to complete the proof of the theorem, we need to show that \eqref{eqn:exchange-poly3} continues to hold with the $L^2(\pa\Omega_k)$-norms replaced by $H^{-1/2}_{\gamma}(\pa\Omega_k)$-norms. Note from Proposition \ref{prop:Steklov-poly} that the eigenvalues, $\lambda^k_j$, of the Dirichlet-to-Neumann operator $T^k_\gamma$ are uniformly bounded below by a constant depending only on $\Omega_k$, and so the $H^{-1/2}_{\gamma}(\pa\Omega_k)$-norm of $Ag$ is bounded from above by $C\gamma^{1/2} \norm{Ag}_{L^2(\pa\Omega_k)}$. Applying Lemma \ref{lem:H-1/2} to $g\in \XM$ therefore implies that
\begin{align*}
    \norm{Ag}_{H^{-1/2}_{\gamma}(\pa\Omega_k)} & \leq C\gamma^{1/2}   \norm{Ag}_{L^2(\pa\Omega_k)}   \leq C\gamma^{1/2} { ( 1 + M^{1/2} )} \gamma^{1/2}\norm{g}_{L^2(\pa\Omega_k)} \\
    & \leq 2C\gamma^{1/2}  ( 1 + M^{1/2} ) \norm{g}_{H^{-1/2}_{\gamma}(\pa\Omega_k)},
\end{align*}
and this completes the proof of the theorem.
\end{proof1}
We are left to prove the bounds on the kernel $K_{\gamma}(x,y)$ from Proposition \ref{prop:kernel-bounds}.
\begin{proof1}{Proposition \ref{prop:kernel-bounds}}
 To prove the kernel bounds, we obtain an explicit representation in terms of a parameterization of the boundary of $\Omega_0$, and then use asymptotics for the modified Bessel function $K_0$. 
 
 Since for any fixed constant $c^*>0$ there is $C>0$ such that
$$|K_0(z)| \leq Ce^{-z}, \qquad z\geq c^*,$$
and likewise for its derivatives, we may assume that the points $x$, $y\in\pa\Omega_0$ satisfy $|x-y|\leq c^*$ for a small constant $c^*$ which will be specified below. Given $x\in\pa\Omega_0$, after a rotation, we may write the boundary of $\pa\Omega_0$ in a neighborhood of $x$ as $y_2 = f(y_1)$. Then, we have, away from any vertex, $n_0 = \frac{1}{\sqrt{1+|f'(x_1)|^2}}(f'(x_1),-1 )$, and so 
\begin{align} \label{eqn:K0}
    K_\gamma(x,y) = 2 \frac{1}{\sqrt{1+|f'(x_1)|^2}}( f'(x_1),-1)\cdot \nabla_x \left(K_0(|x-y|/\gamma)\right).
\end{align}

\emph{Part 1.} In the smooth case, given $x=(x_1,x_2)\in\pa\Omega_0$, after a rotation and translation, we can ensure that $f\in C^4([x_1-c^*,x_1+c^*])$ satisfies,
\begin{align} \label{eqn:f}
 f(x_1) = f'(x_1) = 0, \quad f''(x_1) = \kappa(x),
\end{align} 
in a ball centered at $x$ of radius $c^*>0$. Here, $\kappa(x)$ is the curvature of $\pa\Omega_0$ at $x = (x_1,x_2)=(x_1,0)$. Then, assuming, without loss of generality, that $\Omega_0$ lies above its boundary, the unit outward pointing normal to $\Omega_0$ at $x$ is $-(0,1)$. Therefore, we have
\begin{align*} 
K_\gamma(x,y) = -2\frac{(x_2-y_2)}{\gamma|x-y|}K_0'(|x-y|/\gamma) = 2\frac{y_2}{\gamma|x-y|}K_0'(|x-y|/\gamma).
\end{align*}
Using \eqref{eqn:f}, we can write $y_2 = \tfrac{1}{2}\kappa(x)(y_1-x_1)^2 + O(|y_1-x_1|^3)$, and so
\begin{align} \nonumber
K_\gamma(x,y) & = \frac{\kappa(x)(y_1-x_1)^2 + O(|y_1-x_1|^3)}{\gamma|x-y|}K_0'(|x-y|/\gamma) \\ \label{eqn:K-smooth1}
& = \gamma^{-1}\left(\kappa(x)|x-y| + O(|x-y|^2)\right)K_0'(|x-y|/\gamma).
\end{align}
Here we use the $O$-notation to denote a bound with a constant depending only on $\pa\Omega_0$. The asymptotics of $K_0(z)$ satisfy
\begin{align} \label{eqn:K0-1}
K_0(z) = -\ln(z) + O(1), \quad K_0'(z) = -z^{-1} + O(z|\ln(z)|), \quad K_0''(z) = z^{-2} + O(1)
\end{align}
as $z\to0$, and 
\begin{align} \label{eqn:K0-2}
K_0(z) = \left(\tfrac{\pi}{2z} \right)^{1/2}e^{-z}(1+O(z)),  & \quad K_0'(z) = -\left(\tfrac{\pi}{2z} \right)^{1/2}e^{-z}(1+O(z^{-1})), \\ \nonumber
& K_0''(z) =  \left(\tfrac{\pi}{2z} \right)^{1/2}e^{-z}(1+O(z^{-2}))
\end{align}
as $z\to\infty$. Using these asymptotics for $K_0'(z)$, since $\kappa(x)$ is bounded, from \eqref{eqn:K-smooth1}, for those $y$ with $z=|x-y|/\gamma\leq 1$, the kernel $K_\gamma(x,y)$ is bounded, while for those $y$ with $z=|x-y|/\gamma>1$, the kernel $K_\gamma(x,y)$ satisfies
\begin{align*}
    |K_{\gamma}(x,y)|\leq C\gamma^{-1/2}|x-y|^{1/2}e^{-|x-y|/\gamma} \leq Ce^{-c|x-y|/\gamma},
\end{align*}
as required. In the same way, we can get the required bounds on $\nabla_yK_{\gamma}(x,y)$, by differentiating the equation in \eqref{eqn:K-smooth1}, and using the asymptotics of $K_0'(z)$ and $K_0''(z)$. 

\emph{Part 2.} To handle the polygonal case, we only have to work in a neighborhood of each vertex, $v$, of $\pa\Omega_0$, which after a translation we place at the origin. By the work in the smooth case, we may also assume that $x$ and $y$ lie on different sides of the polygon, and, after a rotation, we assume that $x_1>0$, with $y_1<0$. Moreover, we can choose this rotation such that, in this neighborhood of the vertex, there exists a non-zero constant $\beta$ with $f(y_1) - \beta|y_1|$ equal to a smooth function $F(y_1)$ satisfying $F(y_1) = O(y_1^2) = O(|x_1-y_1|^2)$. Referring back to \eqref{eqn:K0}, and following the method in the smooth case above, we  have
\begin{align*}
   \left| K_\gamma(x,y) - \frac{2}{\sqrt{1+\beta^2}}(\beta,-1)\cdot \nabla_x\left(K_0(|x-y|/\gamma)\right)\right| \leq Ce^{-c|x-y|/\gamma}.
\end{align*}
  Therefore, we are left to show that
\begin{align*}
    (\beta,-1)\cdot \nabla_x\left(K_0(|x-y|/\gamma)\right) & = \frac{1}{\gamma|x-y|}\left( -\beta (y_1-x_1)+(y_2-x_2)) \right)K_0'(|x-y|/\gamma) \\
    &= \frac{1}{\gamma|x-y|}\left( -2\beta y_1+(F(y_1)-F(x_1)) \right)K_0'(|x-y|/\gamma).
\end{align*}
satisfies the estimate in Part 2 of the proposition. Note that we have used $x_2 = \beta x_1+F(x_1)$ and $y_2=-\beta y_1+F(y_1)$ as $x_1>0$ and $y_1<0$. Moreover, since $x$ and $y$ are on different sides of the polygon we have $|x-y|\geq c(|x|+|y|)$. Using the asymptotics of $K_0'$ from \eqref{eqn:K0-1} and \eqref{eqn:K0-2} therefore establishes the required estimate. 
\end{proof1}

\subsection{The estimate in Part 2 of Theorem \ref{thm:exchange} cannot be extended to all of $\HN{s}(\Gamma)$.}\label{sec:Remark1} 
We end this section by showing, as we stated in \eqref{rem:exchange}, that the estimate in Part 2 of Theorem \ref{thm:exchange} cannot be extended to all functions in $\HN{s}(\Gamma)$. For simplicity, we show this in the case that $\pa\Omega_0 = \pa\Omega_1$ is an exact sector in a neighborhood of a vertex placed at the origin.

    Suppose that in a disc of radius $2$, centered at $0$ (which we write as $D_2$), $\Omega_0$ is given by
\begin{align*}
\{(x_1,x_2)\in \R^2: x_2>\beta|x_1|\}
\end{align*}
for a constant $\beta\neq0$. This means that for $\alpha$ satisfying $\tan(\alpha) = \beta^{-1}$, $\Omega_0$ has a vertex of interior angle $2\alpha$. Then, for $x,y\in\pa\Omega_0\cap D_2$, with $x_1>0$,
$$
n_{0}(x)\cdot\nabla (|x-y|) = \frac{-\cos(\alpha)y_1+\sin(\alpha)y_2}{|x-y|} 
= \begin{cases}
    0 &\;\; y_1>0,\\
    \frac{-2\cos(\alpha)y_1}{|x-y|} &\;\;y_1<0.
\end{cases}
$$
Therefore,
$$
(\pa_{n_0,x}\mathcal{G})(x-y)
= \begin{cases}
    0&\;\;y_1>0,\\
    \frac{-2\cos(\alpha)y_1}{\gamma|x-y|}(K_0')\left(|x-y|/\gamma\right) & \;\;y_1<0.
\end{cases}
$$
We will first show that there exists a constant $c>0$, independent of $\gamma$, such that for each $0<\gamma<c$, there exists $f_\gamma\in L^2(\pa\Omega_0)$ such that
\begin{align} \label{examp1}
    \norm{Af_{\gamma}}_{L^2(\pa\Omega_0)} \geq c \norm{f_{\gamma}}_{L^2(\pa\Omega_0)}.
\end{align}
We set $f_{\gamma} = 1_{U_{\gamma}\cap\pa\Omega_0}$, for
\begin{align*}
    U_{\gamma} = \{(x_1,x_2)\in\R^2: x_1\in(-c^*\gamma,0)\}.
\end{align*}
Here $c^*>0$ is a small constant to be chosen below. Then, for $x\in \pa\Omega_0\cap D_2$, with $x_1>0$, we have
\begin{align*}
    Af_\gamma(x) = c_{\beta}\int_{-c^*\gamma}^{0}\frac{y_1}{\gamma|x-y|}(K_0')(|x-y|/\gamma)\,dy_1.
\end{align*}
Here $y = (y_1,-\beta y_1)$ and $c_{\beta}$ is a constant depending only on $\beta$ (that will change from line-to-line). Using the asymptotics of $K_0'(z)$ as $|z|\to0$ from \eqref{eqn:K0-1}, we can write this as
\begin{align} \label{examp2}
      Af_\gamma(x) = c_{\beta}\int_{-c^*\gamma}^{0} \frac{y_1}{ |x-y|^2} \,dy_1+ c_{\beta}\int^0_{-c^*\gamma}\frac{y_1}{\gamma^2|x-y|}O(|x-y|\ln(|x-y|/\gamma))\,dy_1.
\end{align}
For $x=(x_1,\beta x_1)\in\pa\Omega_0$, with $0<x_1<c^*\gamma$, the first integral can be bounded from below by $c_{\beta}\left|\ln(x_1/(c^*\gamma))\right|$, while the second integral can be bounded from above by $C_{\beta}(c^*)^2\left|\ln(x_1/\gamma)\right|$ for a, possibly large, constant $C_{\beta}$ depending only on $\beta$. Since
\begin{align*}
\int_{0}^{c^*\gamma}\left|\ln(x_1/(c^*\gamma))\right|^2 \,dx_1= 2c^*\gamma \qquad \int_{0}^{c^*\gamma} (c^*)^4|\ln(x_1/\gamma)|^2\,dx_1 = O((c^*)^5\gamma),
\end{align*}
for $c^*$ sufficiently small, we have 
\begin{align*}
     \norm{Af_{\gamma}}^2_{L^2(\pa\Omega_0)} \geq c_\beta c^*\gamma = c_\beta\norm{f_{\gamma}}^2_{L^2(\pa\Omega_0)}.
\end{align*}
As $c_\beta$ is independent of $\gamma$, this proves the estimate in \eqref{examp1}.

To complete the construction of  \eqref{rem:exchange}, by interpolation, we just need to show that the lower bound in \eqref{examp1} continues to hold when the $L^2(\pa\Omega_0)$-norms are replaced by the $H^{-1/2}_{\gamma}(\pa\Omega_0)$-norms. Using the lower bound on the Steklov eigenvalues $\lambda_j^k$ from Proposition \ref{prop:Steklov-poly}, the $H^{-1/2}_\gamma(\pa\Omega_0)$-norm of $f_{\gamma}$ is bounded from above by $C\gamma^{1/2}\norm{f_\gamma}_{L^2(\pa\Omega_0)}$. Finally, we use the pointwise lower bound on $Af_{\gamma}(x)$ for $x\in\pa\Omega_0$ with $0<x_1<c^*\gamma$ coming from \eqref{examp2}, and the definition of the $H^{-1/2}_{\gamma}(\pa\Omega_0)$-norm,
\begin{align*}
\norm{p}_{H^{-1/2}_\gamma(\pa\Omega_0)}= \sup\left\{\left|\int_{\pa\Omega_0}pv\right|\,:\,\norm{v}_{H^{1/2}_{\gamma}(\pa\Omega_0)}=1\right\}.
\end{align*}
This provides a lower bound on $\norm{Af_{\gamma}}_{H^{-1/2}_{\gamma}(\pa\Omega_0)}$ comparable to $\gamma$, which gives the required lower bound of 
\begin{align*}    \norm{Af_{\gamma}}_{H^{-1/2}_\gamma(\pa\Omega_0)} \geq c \norm{f_{\gamma}}_{H^{-1/2}_\gamma(\pa\Omega_0)}.
\end{align*}

\section{Application to non-local numerical methods} \label{sec:Claeys}

In this section, we will describe in more detail how the estimates from Theorems  \ref{thm:exchange} and \ref{thm:DtN} relate to the non-local formulation of a Helmholtz scattering problem introduced by Claeys in \cite{Cl}. This formulation, which as we discussed in the introduction can be viewed as a variant of an Optimized Schwarz Method,  is defined for dimension $d\leq 3$, but here we will focus on 2 dimensions.

Claeys studies the wave propagation problem, for $u\in H^1_{loc}(\R^2)$,
\begin{equation}\begin{cases} \label{eqn:wave}
    -\text{div}(\mu\nabla u)-\kappa^2u = f \text{ in } \R^2, \\
    \lim_{\rho\to\infty}\int_{\pa D_{\rho}}|\pa_{\rho}u-i\kappa u|^2\,d\sigma_{\rho} = 0,
\end{cases}
\end{equation}
via a non-overlapping domain decomposition approach. Partitioning $\R^2$ as $\R^2 = \cup_{k=0}^{J}\bar{\Omega}_j$, with $\Omega_k$ bounded for $1\leq k \leq J$, and skeleton $\Gamma = \cup_{k=0}^{J}\pa\Omega_j$, the multi-domain Neumann trace space is defined by
\begin{align*}
    \HN{-1/2}(\Gamma) = H^{-1/2}_\gamma(\pa\Omega_0)\times\cdots \times H^{-1/2}_\gamma(\pa\Omega_J).
\end{align*}

In \cite{Cl} it is shown that this wave propagation problem can be reformulated in terms of non-local exchange and scattering operators:

For $\phi = (\phi_0,\ldots,\phi_J) \in \HN{-1/2}(\Gamma)$, analogously to \eqref{eqn:Psi1}, we write
\begin{align*}
    \PsiN{}(\phi)(x) = \sum_{k=0}^{J}\int_{\pa\Omega_k}\mathcal{G}(x-y)\phi_k(y)\,d\sigma(y),
\end{align*}
and let $\tau_N = (\tau^0_N,\ldots,\tau^J_N)$ be the Neumann traces. In Section 5 of \cite{Cl}, the multi-domain non-local exchange operator $\Pi_\gamma$ is then given by
\begin{align*}
   \Pi_\gamma = Id-2\tau_N\cdot \PsiN{} 
\end{align*}
and is shown to be an isometry on $\HN{-1/2}(\Gamma)$. Here, $\tau_N\cdot \PsiN{}(\phi)$ is the element in $\HN{-1/2}(\Gamma)$ given by evaluating by taking the Neumann trace of $\PsiN{}(\phi)$ from the interior of each subdomain $\Omega_k$. In the case that $J=1$, with $\Omega$ decomposed into two subdomains, this operator $\Pi_\gamma$ is the same as the one defined in \eqref{defn:exchange}.

Letting $T^k_\gamma$ be the DtN operator of $\Omega_k$, for the operator $-\Delta+\gamma^{-2}$, the ingoing and outgoing Robin trace operators are then defined on $H^1_\gamma(\Omega_k)$ by
\begin{align*}
    \tau_{\pm,\gamma}^k (u) = \mu_k\tau_N^k(u) \pm i \omega \gamma^{-1} T^k_\gamma(\tau_D^k(u)).
\end{align*}
Here $\mu_k = \mu|_{\partial \Omega_k}$, $\omega>0$ is an impedance constant, and $\tau_D^k$ is the Dirichlet trace on $\pa\Omega_k$. In \cite{Cl} these traces are used to define a non-local scattering operator $\mathcal{S}_\gamma = \text{diag}_{k=0,\ldots,J}(S_\gamma^k)$, with
\begin{align*}
    S_\gamma^k(\tau^k_{-,\gamma}(u)) = \tau^k_{+,\gamma}(u),
\end{align*}
which is non-expansive on $\HN{-1/2}(\Gamma)$. When $J=1$, $\mu_k=1$ and $\omega = \gamma$, these scattering operators $S_\gamma^k$ are the same as the one defined in Section \ref{subsec:DtN}. These non-local operators $\Pi_\gamma$ and $\mathcal{S}_\gamma$ are then used in Section 7 of \cite{Cl} to reformulate \eqref{eqn:wave} to an equivalent well-posed problem: If $u\in H^1_{loc}(\R^2)$ is the unique solution to \eqref{eqn:wave}, then the tuple of traces $\mathfrak{p}=\tau_{-,\gamma}(u)$ satisfies
\begin{align*}
    \mathfrak{p}-(\Pi_\gamma\cdot \mathcal{S}_\gamma)\mathfrak{p} = \mathfrak{f}, \qquad \mathfrak{p}\in  \HN{-1/2}(\Gamma).
\end{align*}
Here $\mathfrak{f} = \Pi_\gamma(\tau_{+,\gamma}(\phi_f))$, where $\phi_f$ solves the wave-propagation problem in each $\Omega_k$, with zero ingoing Robin data, $\tau_{-,\gamma}(\phi_f) = 0$. As shown in Theorem 7.1 of \cite{Cl} this is a strongly coercive formulation of \eqref{eqn:wave}.
This formulation is valid for each fixed $\gamma>0$, and the operators $\Pi_\gamma$ and $\mathcal{S}_\gamma$ are non-local for each $\gamma$.

Our estimates from Theorems \ref{thm:exchange} and \ref{thm:DtN} then give information on the extent that these operators resemble local operators in the $\gamma\to0$ limit. If, for example, each piece of the partition $\Omega_k$ is contained in the interior of $\Omega_{k-1}$, then we can be in the setting where the boundaries of each subdomain are smooth. However, if the partition has any cross-points (where three or more of the subdomains meet at a point), then necessarily some of the boundaries of the subdomains will be polygonal, with convex vertices.

As a side remark, we note that the scattering operator $\mathcal{S}_0$ bears similarity to the impedance-based merge operators that arise in the hierarchical Poincaré–Steklov method. In that framework, the domain is decomposed into a hierarchical tree of subdomains, and a merge operator is employed to reconstruct the global solution. A version of this method was introduced by Gillman, Barnett, and Martinsson \cite{gillman2015spectrally} to numerically approximate solutions to Helmholtz-type equations in the context of potential scattering. It has also been applied to obstacle scattering problems \cite{pedneault2017schur}. The merge operator was further studied in \cite{beck2022quantitative}, where function spaces analogous to the spaces $X_\gamma(M)$ introduced here were defined in terms of the scattering frequency, and were used to prove boundedness of the merge operator.

We first consider the scattering operator $\mathcal{S}_\gamma = \text{diag}_{k=0,\ldots,J}(S_\gamma^k)$. Defining the local version of this operator via
\begin{align*}
    S^k_{0}(\tau_{-,0}^k(u)) = \tau_{+,0}^k(u), \qquad  \tau_{\pm,0}^k(u) = \mu_k\tau_N^k(u) \pm i \omega \gamma^{-1}\tau_D^k(u)
\end{align*}
we have
\begin{align*}
  S^k_{0}(\tau_{-,0}^k(u))  -  S^k_\gamma(\tau_{-,\gamma}^k(u)) = \tau_{+,0}^k(u) - \tau_{+,\gamma}^k(u) = i \omega \gamma^{-1}\tau_D^k(u) - i \omega \gamma^{-1}T_\gamma^k(\tau_D^k(u)).
\end{align*}
Therefore, if the Dirichlet trace $\tau_D^k(u)$ is in $\XM$, then as in Corollary \ref{cor:DtN}, we can bound this expression. For example, from Theorem \ref{thm:DtN}, for each $0<a<1$ there exists a constant $C_a$ such that for $M <C_a^{-1}\gamma^{-a}$,
    \begin{align*}
        \norm{S^k_{0}(\tau_{-,0}^k(u))  -  S_\gamma^k(\tau_{-,\gamma}^k(u))}_{H^{-1/2}_{\gamma}(\pa\Omega_k)} \leq C_a\omega\gamma^{-1}(\gamma^{1/3} + M^{1/2}\gamma^{a/2})\norm{\tau_D^k(u)}_{H^{-1/2}_{\gamma}(\pa\Omega_k)}.
    \end{align*} 
In particular, for an impedance constant $\omega$ comparable to $\gamma$, the right hand side is $o(1)$ as $\gamma\to0$. If the sub-domains in the partition have smooth boundaries (so that in particular, there are no cross-points), then one can improve this estimate using the no convex vertex and smooth case estimates from Theorem \ref{thm:DtN}.

We now turn to the non-local exchange operator $\Pi_\gamma$. Denoting $\Gamma_{k\ell}$ to be the (possibly empty) shared boundary between $\Omega_k$ and $\Omega_\ell$, we define the operator $S^{*,j}_{k\ell}:H^{-1/2}_{\gamma}(\Gamma_{k\ell})\to H^{-1/2}_{\gamma}(\pa\Omega_j)$ by
\begin{align*}
    S^{*,j}_{k\ell}f(x) = \int_{\Gamma_{k\ell}}(\pa_{n_j,x}\mathcal{G})(x-y)f(y)\,d\sigma(y).
\end{align*}
Note that using the kernel bounds from Proposition \ref{prop:kernel-bounds}, as in the proof of Theorem \ref{thm:exchange}, for each $0<a<1$, there exists a constant $C_a>0$ such that for all $0<\gamma<c$, $M<C_a^{-1}\gamma^{-a}$, and $f\in \XM$,
\begin{align} \label{eqn:claeys1}
    \norm{ S^{*,j}_{k\ell}f}_{H^{-1/2}_\gamma(\pa\Omega_j)} \leq C{(1+M^{1/2})}\gamma^{1/2}\norm{f}_{{H}^{-1/2}_\gamma(\Gamma_{k\ell})}.
\end{align}
Again, this estimate can be improved and on the whole of $H^{-1/2}_\gamma(\Gamma_{k\ell})$ if the boundaries are smooth. Following the analogous calculation to the proof of Lemma \ref{lem:Pi} above, for example, for $x\in \Gamma_{01}$, the zero-th component of $\Pi(\phi)$ is given by
\begin{align*}
  -\phi_1(x) - 2\sum_{k=0}^{M} \sum_{\ell\neq k} S^{*,0}_{k\ell}\phi_k(x).
\end{align*}
Analogously to \eqref{eqn:Pi0}, we denote $\Pi_0$ to be a local exchange operator, exchanging Neumann data on each $\Gamma_{jk}$ up to a change in sign. Then, from the estimates in \eqref{eqn:claeys1}, there exists a constant $C_a>0$ such that for all $0<\gamma<c$, $M<C_a^{-1}\gamma^{-a}$, and for $\phi = (\phi_0,\ldots,\phi_J)\in \HN{-1/2}(\Gamma)$ with $\phi_k\in \XM$, we have
\begin{align*}
    \norm{(\Pi_\gamma-\Pi_0)\phi}_{\HN{-1/2}(\Gamma)} \leq C{(1+M^{1/2})}\gamma^{1/2}\norm{\phi}_{\HN{-1/2}(\Gamma)}.
\end{align*}
Again the right hand side is $o(1)$ as $\gamma\to0$, and if the subdomains have smooth boundaries, then one obtains an improved estimate on the whole of $\HN{-1/2}(\Gamma)$.

\bibliographystyle{plain}
\bibliography{nonlocal}

\end{document}